\documentclass[12pt]{amsart}
\usepackage{epsfig}

 \newcommand{\new}{\newcommand}                        
 \new{\trunc}{\hat{\otimes}}                           
 \new{\tensor}{\otimes}                                
 \new{\iso}{\cong}                                     
 \new{\W}{\mathfrak{W}}                                
 \new{\g}{\mathfrak{g}}                                
 \new{\uqg}{U_q(\g)}                                   
 \new{\bracket}[1]{\langle#1\rangle}                   
 \new{\qdim}{\operatorname{qdim}}                      
 \new{\Hom}{\operatorname{Hom}}                        
 \new{\Ad}{\operatorname{Ad}}                          
 \new{\id}{{\operatorname{id}}}                        
 \new{\tr}{{\operatorname{tr}}}                        
 \new{\C}{\mathbb{C}}                                  
 \new{\Cat}{\mathfrak{C}}                              
 \new{\Z}{\mathbb{Z}}                                  
 \new{\R}{\mathbb{R}}                                  
 \new{\Q}{\mathbb{Q}}                                  
 \new{\defequals}{\stackrel{\rm def}{=}}               
 \newcounter{letter}                                   
 \newenvironment{alist}{
 \begin{list}{{(\alph{letter})}}{\usecounter{letter}}
 }{\end{list}}                                         
 \newtheorem{thm}{Theorem}                             
 \newtheorem{prop}{Proposition}                        
 \newtheorem{cor}{Corollary}                           
 \newtheorem{lem}{Lemma}                               
 \newtheorem{defn}{Definition}                         
\theoremstyle{definition}                              
\newtheorem{rem}{Remark}                               
\new{\nonsimply}{Sawin00a}                             
\new{\framed}{Sawin00d}                                
\new{\closed}{Sawin00c}                                
\new{\pic}[5]{\raisebox{#3pt}{
\hspace{#4pt}\epsfig{file=#1.ps,height=#2pt}\hspace{#5pt}}}

\begin{document}
\title[Invariants of Spin Three-Manifolds]{Invariants of Spin
Three-Manifolds From {Chern-Simons} Theory 
  and  Finite-Dimensional Hopf Algebras}
\author{Stephen F. Sawin}
\address{Fairfield University\\(203)254-4000x2573\\
ssawin@fair1.fairfield.edu}

\begin{abstract}
A version of Kirby calculus for spin and framed three-manifolds
 is given and is
used to construct invariants of spin and framed three-manifolds in two
situations.  The first is ribbon $*$-categories which possess odd
degenerate objects.  This case includes the quantum group situations corresponding
to the half-integer level Chern-Simons theories conjectured to give
spin TQFTs by Dijkgraaf and Witten \cite{DW90}.   In particular, the
spin invariants constructed by Kirby and Melvin \cite{KM91} are shown
to be identical to the invariants associated to $\mathrm{SO}(3).$ 
Second, an invariant of spin  manifolds analogous to the Hennings
invariant is constructed beginning with an arbitrary factorizable, unimodular
quasitriangular Hopf algebra.  In particular a framed manifold invariant is
associated to every finite-dimensional Hopf algebra via its quantum
double, and is conjectured to be identical to Kuperberg's
noninvolutory invariant of framed manifolds associated to that Hopf algebra.
\end{abstract}
\maketitle

\section*{Introduction}

This article is motivated by, and addresses, two separate questions.
The first is Dijkgraaf and Witten's remarkable observation in
\cite{DW90}  based on the path-integral formulation.  They argue that for certain
nonsimply-connected Lie groups $G$ the level $k$ in the definition of
the Chern-Simons field theory,  which ordinarily must be an
integer in order to get a well-defined topological quantum field
theory, can be a half-integer and still be expected to yield a
sensible theory in the spin category.  One would like to be able to
reproduce this observation rigorously in the algebraic/combinatorial
quantum group formulation, as part of a general effort to relate
these two settings. 

 The second question has to do with the mysterious
`nonsemisimple' topological invariants coming from quantum groups
discovered by Kuperberg \cite{Kuperberg96} (generalizing the
semisimple version studied by many authors including Kuperberg
\cite{Kuperberg91}, Barrett and Westbury \cite{BW95,BW96}, and Chung, Fukuma, and
Shapere
\cite{CFS94})   and Hennings 
\cite{Hennings96}  (also considered by Lyu\-ba\-shen\-ko \cite{Lyubashenko95b,Lyubashenko95d},
Kerler \cite{Kerler95,Kerler97}, and Kauffman and Radford
\cite{Kauffman95,KR95b}), which appear to be closely related to, but distinctly
different from, the Chern-Simons invariants.  The Kuperberg invariant,
which assigns an invariant of framed three-manifolds to each
finite-dimensional Hopf algebra $H,$ is widely conjectured to be
identical to the Hennings invariant associated to the quantum double
of $H.$
Unfortunately, for a typical Hopf algebra the quantum double is
quasitriangular, but not necessarily ribbon, while the Hennings
invariant requires a ribbon Hopf algebra (with some additional
nondegeneracy conditions).  This additional structure is reflected in
the fact that the Hennings invariant depends only on a $2$-framing
(which can be  normalized away with some loss in terms of cutting
and pasting relations), rather than  the more involved framing.  Thus one
would like to extend the Hennings invariant to associate
framing-dependent invariants to quasitriangular Hopf algebras.   This is 
especially important because the framing structure represents much of
the complexity that makes the Kuperberg invariant difficult to work
with.

This article solves both of these questions in a common framework.
Relying on the fact that, crudely speaking, a framing is a spin
structure plus a $2$-framing, we identify in both the Chern-Simons and 
nonsemisimple case the weakening of
algebraic structure from ribbon to quasitriangular with a weakening of
the topological invariance from dependent on a $2$-framing to dependent
on a framing.  This identification can be traced ultimately to the
link invariant level, where in the quasitriangular case the natural
quantities that arise are invariants only of links with even framing
or self-linking number.  

As interesting as the two questions are separately, the connection
between them revealed by this common framework deserves attention
also.  The hints of the geometry of Chern-Simons theory which pervades
the Kuperberg and Hennings invariants seem to demand a physical
explanation, and it is to be hoped that a link between Kuperberg and
Hennings'  algebra on the one hand and Dijkgraaf and Witten's geometry 
on the other offers
a useful step towards such an explanation.

Section 1 gives a combinatorial description of framed three-manifolds
(i.e., equipped with a spin structure and an even $2$-framing)  in
terms of surgery on links all of whose components have even framings.
Section 2 gives a general framework analogous to Reshetikhin and
Turaev's modular Hopf algebras \cite{RT91} for constructing invariants
of spin manifolds and identifies the invariants of this sort arising
from quantum groups with the Chern-Simons theories meeting Dijkgraaf
and Witten's conditions for spin theories.  Also, the $\mathrm{
SO}(3)$ theory is
identified with Kirby and Melvin's \cite{KM91} spin invariants constructed
from quantum $\mathrm{ su}_2,$ and a formulation generalizing theirs is given
in general.  Finally Section 3 defines a Hennings-type invariant of
spin (or framed) three-manifolds starting from a factorizable
quasitriangular Hopf algebra (these conditions include the quantum
double as a special case).  Below we offer a more detailed introduction to each of the
three sections.

\subsection*{Even links, spin three-manifolds and surgery}

Our approach to topological invariants is the now familiar one of what
is sometimes called 
quantum topology:  Describe the topological entity by some
combinatorial data modulo relations generated by a few simple `moves,'
translate the data to algebraic objects, and observe that, magically,
the moves are algebraic relations which the objects at hand satisfy.  Of
course the magic reflects a hidden and poorly understood geometric
underpinning, which we will discuss in the next section.

The first instance of this approach is even links: links with an even
framing on each component.  Here the data is a link projection with
each component having winding number one, and the moves are Kauffman's
regular isotopy, weaker than the usual Reidemeister or framed Reidemeister moves.  In
the  translation of these moves to the language of Hopf
algebras we will see that the ribbon conditions are no longer
necessary.  The use of winding number one rather than zero corresponds
on the algebraic side to evaluating with quantum characters rather
than tracial functions, which is fairly natural from the point of view
of Hopf algebras themselves but does not 
fit neatly in the language of rigid braided categories.

Surgery on framed links produces four-manifolds with three-manifold
boundary.  With the appropriate moves ($2$-framed Kirby moves) 
surgery gives a
combinatorial description of three-manifolds equipped with the
$2$-framing needed to regularize Chern-Simons theory (see Atiyah
\cite{Atiyah90}, and \cite{\framed}).  Surgery on even 
links produce spin four-manifolds with spin three-manifold boundary.
With  the appropriate moves (spin Kirby moves) such an even surgery give a combinatorial
description of spin three-manifolds equipped with a compatible
$2$-framing.  Specifically, the $2$-framing has to equal Rohlin's
$\mu$ invariant modulo $16.$  Since framings of three-manifolds are in
one-to-one correspondence with spin structures together with
$2$-framings, we get a combinatorial description of framed manifolds,
but the compatibility puts a constraint on the possible framings.  The
exact significance of this restriction as it relates to the framings
of Kuperberg's invariant is not clear, although since a fixed shift in
the $2$-framing multiplies each invariant by a fixed quantity, the
distinction is fairly minor.  This section also offers a Fenn-Rourke
version of the spin Kirby moves which  replaces the general handle-slides
with a more restrictive `semilocal' move.

\subsection*{Reshetikhin-Turaev type spin invariants}

Ribbon categories (see Reshetikhin and Turaev \cite{RT91}, Turaev
\cite{Turaev94}, or Kassel \cite{Kassel95}) are the appropriate place
to look for framed link invariants, and if they meet some additional
conditions (summarized in the definition of a modular category in
\cite{Turaev94}) they give invariants of $2$-framed  three-manifolds
through the surgery presentation.  One might think that modularity is
a very restrictive constraint on a ribbon category, but in fact it is
not.  M\"uger \cite{Muger99} and Brugui\`eres \cite{Bruguieres99} give
a kind of quotient of a ribbon $*$-category  (Brugui\`eres works with a
more general  but perhaps less natural substitute for the $*$-structure) that
roughly speaking deletes the part of the category which provides no link
information.  Sometimes this quotient results in a modular category.  We say
sometimes because the objects which stand in the way of modularity (called degenerate)
come in two flavors, odd and even, and only the even can be
eliminated.  Thus the existence of odd degenerate objects is the only
obstruction to constructing a modular category from a ribbon category.

The quotient fails to go through for odd degenerate objects because
they actually carry some information about the link.  In fact they
contribute a factor of $-1$ raised to the self-linking number for every
component they label.  Of course if we restrict our attention to even
links, there is no information at all, and the quotient can go
through.  Imitation of the usual Reshetikhin-Turaev construction in
this case yields an invariant of spin and compatibly-framed three-manifolds (presumably it
yields an appropriately modified version of the axioms of topological
quantum field theory, but we defer that important question  and
focus only on the invariants in this article).  In \cite{\nonsimply} 
the author analyzed the quotient construction applied to subsets of 
the Weyl alcove.   This analysis  allows us to identify the levels at 
which we get a spin
Chern-Simons theory associated to a given simple Lie group.

	On the geometric side, we generalize Dijkgraaf and Witten's
observation, relying on the integration of the generating class of
$H^4(BG,\Z)$ on spin four-manifolds, to a classification of when 
the physical interpretation leads us to  expect  spin theories.  In fact
we get  complete agreement with the algebraic answer (there is
actually a small subtlety:  as in \cite{\nonsimply},  the
Dijkgraaf-Witten theories often factor as a  product of invariants,
but the set of all these factors 
is a complete list of theories arising from quantum groups).

Finally, we imitate Kirby and Melvin's construction of spin invariants
from quantum $\mathrm{ su}_2$ at certain levels to arbitrary quantum
groups.  In fact we  find that it works in just the situations in which there
are Dijkgraaf-Witten theories, and prove that the Kirby-Melvin spin
invariants in fact coincide with the Dijkgraaf-Witten theories.
Kirby and Melvin did their computation in quantum $\mathrm{su}_2,$  and thus
associated their invariant with $\mathrm{SU}(2),$  while we do the
computation entirely in the set of representations associated to
$\mathrm{SO}(3),$ and thus find the invariant more naturally
associated to  $\mathrm{SO}(3)$ as expected  
from the geometry and physics.

\subsection*{Hennings-type invariants of spin manifolds}

In the combinatorial version of Chern-Simons and the other traditional
quantum invariants, one does not work with the quantum group itself, which
is not semisimple, nor with the whole of the representation theory.
Instead one works  with a piece of the representation theory modified so as to
resemble the representation theory of a semisimple ribbon Hopf
algebra.  The theory outside the Weyl alcove, which roughly
corresponds to the nonsemisimple part of the quantum group, is simply
thrown away.  

The Hennings  invariant, by contrast, relies heavily on
the nonsemisimple part of the quantum group.   More precisely,
 the link invariants naturally
associated to the three-manifold invariant are labeled by quantum
characters, which include the quantum traces of irreducible
representations  labeling links in Reshetikhin and Turaev approach,
but also include functions associated to nonirreducible
representations.  What's more the analogue of the surgery label 
(called in the sequel  $\omega$) which Reshetikhin and Turaev use to construct the
three-manifold invariant is the left integral.  This integral  turns out to be
an element of the socle of the algebra of quantum characters  
(roughly,  a maximally nilpotent element).   The topological effect
of this nonsemisimplicity is that the invariant is zero except for
rational homology three-spheres and satisfies only a subset of the
cut-and-paste axioms of TQFT which its semisimple cousins satisfy.

We find for a quasitriangular Hopf
algebra which is not necessarily ribbon that quantum characters label invariants  
of even links  (In the
presence of a ribbon element, quantum characters and cocommutative
functionals can be used almost interchangeably, but in our more
delicate situation we find that quantum characters play the more
fundamental role).  If the algebra is unimodular, the left integral in
the dual (which is also the right integral) is a quantum character,
and in the presence of a nondegeneracy condition (factorizability
suffices) it has the appropriate properties to give an invariant of
compatibly framed three-manifolds.  The quantum double of a finite-dimensional
Hopf algebra meets both conditions (factorizable and unimodular),
though it is often not ribbon, and we conjecture that the invariant of
framed manifolds which we associate to the quantum double of $H$ is
the framed invariant Kuperberg associates by very different means to
$H.$  We also conjecture that the Dijkgraaf and Witten spin invariants
of the previous section arise by a construction analogous to that of
Reshetikhin and Turaev from a quasitriangular but not ribbon Hopf
algebra.

\subsection*{Acknowledgements}  I would like to thank J. Baez, J.
Barrett, R. Kirby, A. Liakhovskaia, G. Masbaum,  E. Witten and G.
Zuckermann for helpful conversations and suggestions.

\section{Even Links, Spin Three-Manifolds and Surgery}
\subsection{Even Links}

A framed link is an oriented link in oriented $S^3$ together with a nonzero section of
the normal bundle, considered up to ambient isotopy.  The framing or
self-linking number of a component is the linking number between that
component and a copy of it pushed off slightly in the direction of the
framing.  An \emph{even} link is a framed link all of whose components have
even self-linking number.  A projection of a framed link (in
particular of an even framed link) is a projection of any
representative of the isotopy class onto the oriented plane  (the
representative is first viewed as sitting inside $\R^3$) such that
the projection of the link is  a smooth immersion of the union of
circles with no self-intersections other than transverse double
points, and the
framing is never orthogonal to the projection.
Of course such a projection, together with identification of each
crossing as over or under, uniquely determines the isomorphism class
of the framed link.    We will sometimes  be interested in such projections
which also come equipped with a \emph{height function,}  i.e., a
smooth map from $\R^2$ to $\R$ without critical points such that
when the map is restricted to the projection  its critical points  are
nondegenerate and do not fall on the crossings.   

An isotopy of a framed link and  a projection (and of a height
function) is collectively called a \emph{planar isotopy}
(respectively \emph{simple isotopy}) if at each
point of the isotopy the link projection (and height function) satisfy
the conditions of the previous paragraph.

The winding number of a component of a projection is the total number
of complete counterclockwise rotations the tangent vector to that
component of the 
projection undergoes in a complete circuit around the component in the
direction of its orientation.   It
is also half the total signed number of critical points of the height
function, with those turning counterclockwise counting as $+1$ and
those turning clockwise counting as $-1.$  

Note that the winding number of a component of a projection and the
framing of that projection are always of opposite parity  (i.e. one is
odd, the other is even).   To see this notice each changes parity only
when the other does under the Reidemeister moves (see Burde and
Zieschang \cite{BZ86}), so
if the claim is true for one projection of a link it is true for all
projections of that link.   It is certainly true for a component which 
is a framed  unknot
unlinked with the other components (using the simple projection in
which it participates in no crossings with other components and 
crosses with itself $n$ times where $n$ is the absolute value of the 
framing) and it remains true
when a single crossing in a projection is switched from over to under
or under to over.  Since it is well known (see Adams \cite{Adams94}, for 
example)  that one can get from any projection to a projection of
unlinked unknots by a sequence of such crossing changes, the claim is
true in all projections.

\begin{figure}[hbt]
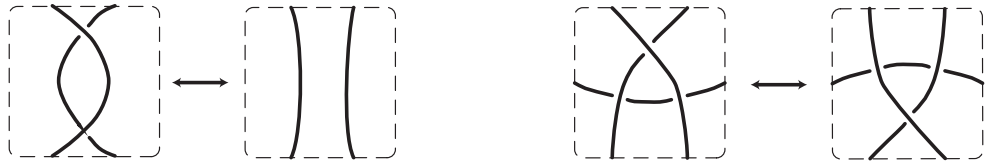

$$\pic{regular}{60}{-25}{0}{0}$$
\caption{The moves of regular isotopy}\label{fg:Reidemeister}
\end{figure}

\begin{figure}[hbt]
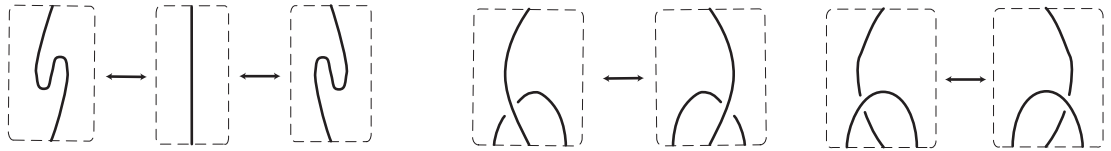

$$\pic{height}{55}{-25}{0}{0}$$
\caption{Moves which change the height function}\label{fg:height}
\end{figure}

\begin{prop}\label{pr:Reidemeister}
Every link admits a projection in which each  component with even framing has
winding number  one and each component with odd framing has winding
number zero.   Two such projections are of the same framed link if and
only if they can be connected by 
planar isotopy together with the \emph{regular isotopy  moves} shown in Figure
\ref{fg:Reidemeister} (understood to apply with any orientations on
the pictured
strands).   Two such projections  equipped with height
functions can be connected by simple isotopy together with the regular
isotopy  moves (with
the vertical indicating the height function) and  the \emph{height
function moves} in Figure
\ref{fg:height}.
\end{prop}

\begin{figure}[hbt]
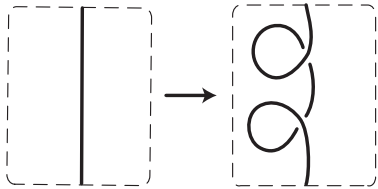

$$\pic{doubletwist}{70}{-30}{0}{0}$$
\caption{Changing the winding number without changing the framing
}\label{fg:doubletwist} 
\end{figure}

\begin{proof}
Recall every framed link admits a projection (every link admits a
projection according to \cite{BZ86}, and by adding a suitable number
of full twists to the projection one can make the projection have any
given framing), so choose one such.  Notice the moves
in Figure \ref{fg:doubletwist} do not change the framed link, but
change the winding number of the indicated component by $\pm 2.$  A
sequence of such moves gives a projection in which every component has
winding number one or zero, and the observation preceding the
proposition indicates these components must have respectively even or
odd framings.

For the second claim, Trace \cite{Trace83} has shown that two link
projections are regular isotopic  (i.e., connected by a sequence of
regular isotopy moves and planar isotopy), if and only if their links
are isotopic and each component has the same writhe and winding number.  The
writhe of a component of a projection is the signed sum of the
crossings of the component with itself, which is of course the framing
if it is a projection of a framed link.  

That every framed link admits a height function is clear, by choosing
any real function with no critical points and perturbing slightly as
necessary.  That two such are connected by a sequence of the height
function moves appears in Kassel \cite[Lemma X.3.5]{Kassel95}.
\end{proof}

\begin{rem}
The appropriate moves for framed links without the restriction on winding
number are the regular isotopy moves of Figure \ref{fg:Reidemeister} together
with a move that replaces a positive full twist with winding number $+1$ with a
positive full twist with winding number $-1.$  These are called the framed
Reidemeister moves \cite{Kassel95,Sawin96a}.  The regular isotopy
moves  are strictly weaker.
\end{rem}

In Section 3 we will assign numbers to link projections with each component
having winding number one, such that the numbers are unchanged by simple
isotopy and the moves of regular isotopy, and thus are invariants of the even
links the projections represent.

\subsection{Framed and Spin Three-Manifolds}

Recall that the ordinary Chern-Simons invariant depends not simply on a
three-manifold but on a $2$-framed three-manifold.  A $2$-framing
\cite{Atiyah90,\framed} of a three-manifold $M$ can either be defined 
as a  trivialization up 
to isotopy
of a $\mathrm{Spin}(6)$ bundle associated  with $TM \oplus TM,$ or as 
a choice of cobordism class of four-manifolds which $M$ bounds 
\cite{Kerler99}.   We will use the second definition, and recall that 
the cobordism class of a four-manifold is determined by its boundary 
and signature, so we may think of a $2$-framing on $M$ as an integer 
representing the signature of the four-manifold \cite{Kerler99}.

The choice of a spin structure  (recall every compact oriented
three-manifold admits a spin structure,  defined as 
a lifting  up to isotopy of the $\mathrm{SO}(3)$ bundle associated to
$TM$ to an $\mathrm{SU}(2)$ bundle) and a  $2$-framing for $M$ is
equivalent to the choice of a framing for $M$, i.e. a trivialization
of the tangent bundle up to isotopy \cite{\framed}.    That is, if two 
framings induce the same spin structure and $2$-framing they are 
homotopic, and every combination of spin structure and $2$-framing is 
induced by some framing.

According to Kirby \cite[Chapters II, IV, and VII]{Kirby89} every spin three-manifold spin bounds 
a spin  four-manifold, in
fact a $2$-handlebody (i.e., a four-manifold formed by attaching a 
collection of $2$-handles to a $0$-handle).  A $2$-handlebody is spin 
if and only if its intersection form on second cohomology is even, 
and in this case possesses a unique spin structure.  Thus we may 
specify a spin structure on a three-manifold by specifying a 
$2$-handlebody with even form which the three-manifold bounds.   Now a $2$-framing on a three-manifold
can also be determined by specifying a four-manifold which it bounds,
in fact it is determined by the signature of the intersection form 
(again,  it may as well be a $2$-handlebody).  It
is natural to represent both pieces of information by a single
four-manifold, but this is only possible for certain framings.  
Recall \cite[Chapter XI]{Kirby89}
that if $M$ is a spin three-manifold then Rohlin's invariant $\mu(M)$
is an even integer modulo $16$ such that every four-manifold which
spin bounds $M$ has signature equal to $\mu(M)$ modulo $16.$
Motivated by this we define

\begin{defn}
A framing on a compact, connected, oriented three-man\-i\-fold $M$ is
called \emph{compatible} if the induced $2$-framing as an integer is equal
to $\mu(M) \mod 16.$
\end{defn}

Thus a compatible framing on $M$ can be represented by a  
$2$-handlebody.

\subsection{Spin Manifolds and Surgery}

Recall from Kirby \cite[Chapter I]{Kirby89} that if $W$ is a 
$2$-handlebody, 
 $W$ can be presented by an unoriented framed link in $S^3.$  Here $S^3$ bounds the
$0$-handle $B^4,$ and the link with each component thickened to a
solid torus with distinguished longitude determines how to attach a
$2$-handle along each component.  The matrix of the intersection form
is given by the linking matrix, which is even (i.e. has even entries
along the diagonal) if and only if the $2$-handlebody is spin.  The
boundary of a four-manifold described by a link is
the three-manifold obtained by surgery on that link.    That is, it is the 
result of  removing
a tubular neighborhood of each component and gluing it back in by sending
the meridian to the longitude and the longitude to minus the
meridian.  Thus surgery on an even link can be viewed as resulting in
not simply a three-manifold, but a three-manifold together with a
compatible  framing.

The following theorem and proof are direct translations of Kirby's
well-known surgery theorem \cite{Kirby78} to the spin case.

\begin{figure}[hbt]
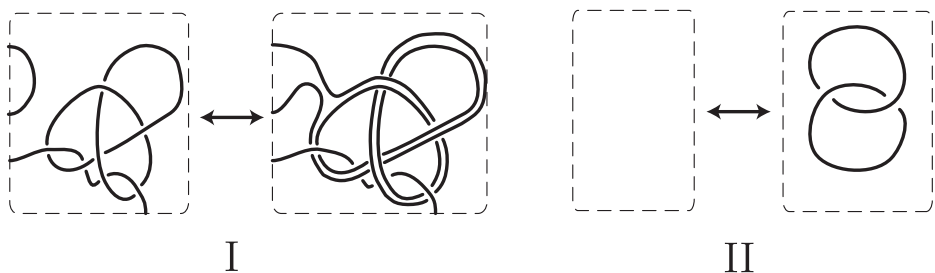

$$\pic{kirby}{100}{-45}{0}{0}$$
\caption{Spin Kirby moves}\label{fg:Kirby}
\end{figure}

\begin{thm}\label{th:Kirby}  \mbox{}
\begin{alist}
\item
Every three-manifold with a compatible framing can be presented by surgery
on an even link.  Two such presentations determine the same framed
three-manifold if and only if they can be connected by a sequence of
the spin Kirby moves pictured in Figure \ref{fg:Kirby}.  Here move I
is the usual handleslide or band connect sum of Kirby's original
theorem.
\item  Every spin three-manifold can be presented by surgery on an even
link, and two such determine the same spin three-manifold if and only
if they can be connected by a sequence of spin Kirby moves I and II
as above together with distant union with the link representing the Kummer
surface, pictured for example in \cite[page 9]{Kirby89}.
\end{alist}
\end{thm}

\begin{proof} \mbox{}
\begin{alist}
\item  The first sentence is clear, as is the fact the the spin
Kirby moves do not change the framed three-manifolds.  So suppose
$M^4$ and $N^4$ are two spin $2$-handlebodies with the same spin boundary and
signature presented by links $L_M$ and $L_N.$  We will show that
$L_M$ and $L_N$ can be connected by spin Kirby moves.

Since they each have the same signature,  gluing $M$ to $-N$ along
$\partial M \times I$ we get a closed spin manifold with signature
zero.   Hence by \cite[Thm. VII.3]{Kirby89} there is a spin five-manifold $W$
which this closed manifold bounds.  Choose a Morse function $f:W\to [1,2]$
such that 
\begin{align*}
f^{-1}[1]&=M\\
f^{-1}[2]&=N\\
f(\sigma,t)&=t \,\,\text{ for } \,\,(\sigma,t) \in \partial M \times I.
\end{align*}

Each $f^{-1}[t]$ for $t$ not a critical point is a spin
four-manifold.  If $t$ is an index $1$ critical point then
$f^{-1}[t+\epsilon]$ is $f^{-1}[t-\epsilon]$ connect summed with $S^1
\times S^3,$ with the same spin structure outside the region of the
connect sum.  In between these two manifolds we can change $W$ by
replacing the $1$-handle with a $2$-handle attached along a
contractible loop in the boundary.  The spin structure on the boundary
$f^{-1}[t \pm \epsilon]$ induces  a spin structure on the boundary
$S^4$ of the $2$-handle $B^5,$ which extends to the interior.    In
this fashion $W$ can be replaced by a new spin five-manifold in which
there are no $1$-handles but the Morse manifolds  (i.e. $f^{-1}$ of
noncritical values) are the same.  Similarly we can replace all the
$4$-handles with $3$-handles.

Arrange $f$ so that all $2$-handles have Morse values less than those
of all $3$-handles.   Because the Morse manifolds  are simply-connected,
each $2$-handle as it is attached connect sums $S^2 \times S^2$ to the Morse manifold,
which corresponds to spin Kirby move II applied to $L_M.$  Likewise by 
flipping $f$ we see 
each $3$-handle corresponds to move II applied to $L_N.$  Thus a sequence of Moves II
applied to $L_M$ and $L_N$ gives links representing the same spin
four-manifold with boundary $\partial M.$  Kirby's argument that these
can be connected by a sequence of moves I goes through unchanged.
\item If $M$ and $N$ have the same spin three-manifold boundary, then their
signature differs by a multiple of $16.$  Thus the union of one of
their links
with sufficiently many copies of the Kummer surface link results in
two links which present the same framed manifold, and part (a) of the
theorem  applies. \qed
\end{alist}
\renewcommand{\qed}{}
\end{proof}

\begin{rem}
The link invariants we construct will be multiplicative in the sense
that the invariant of a distant union of links  (i.e.  the union of two links
embedded simultaneously in $S^3$ so that they are separated by an
$S^2$) is the product of their individual invariants.  If a
multiplicative link invariant $I(L)$ is invariant under the spin Kirby
moves, it is an invariant of compatibly framed three-manifolds.  
Furthermore $I(L) I(K)^{-\sigma/16},$ where $K$ is an even
link representing the Kummer surface and $\sigma$ is the signature of
the linking matrix of $L,$ is invariant under the additional move of
part (b) and thus gives an invariant of spin manifolds.  Since the
spin and framed versions of the invariant differ only by this simple
normalization, we will move freely between the two versions.  The two
versions are exactly analogous to the ordinary and $2$-framed version
of the three-manifold invariant of Reshetikhin and Turaev
\cite{RT91}.  As in that case,  we expect the framed version to be
more natural at the level of TQFTs.
\end{rem}

Just as Fenn and Rourke's \cite{FR79} `semilocal' simplification of
Kirby's surgery theorem gives an alternate set of moves which are
sometimes more convenient for addressing three-manifold invariance,
so we will find it helpful to have the following version of the spin
Kirby moves at hand.

\begin{figure}[hbt]
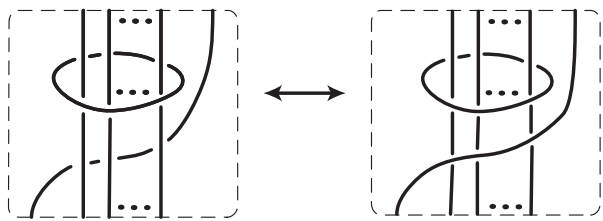

$$\pic{spinFR}{80}{-35}{0}{0}$$
\caption{Alternate spin Kirby move I$'$}\label{fg:spinFR}
\end{figure}

\begin{prop}\label{pr:spin}
Spin Kirby moves I and II of Figure \ref{fg:Kirby} generate the same
equivalence relation on links as moves I$\,'$ and II, where I$\,'$ is pictured
in Figure \ref{fg:spinFR}, with  the number of strands passing through
the unknot being arbitrary.
\end{prop}

\begin{figure}[hbt]
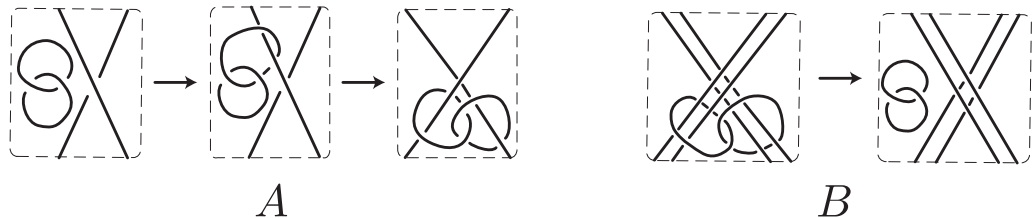

$$\pic{switch}{80}{-35}{0}{0}$$
\caption{Using handle-slides with the Hopf link to switch a crossing}
\label{fg:crossingswitch} 
\end{figure}

\begin{proof}
Of course it suffices to take an arbitrary instance of spin Kirby move
I and decompose it as a sequence of moves I$'$ and II.

Let $L$ be an even  link with a component $C$ and let $B$ be another
component to be band connect summed with $C,$ as illustrated in Kirby 
move I in
Figure \ref{fg:Kirby}.  Choose a presentation of $L$ with the winding
number of $C$ being one as in Proposition \ref{pr:Reidemeister}, and such that the
band between $B$ and $C$ along which the connect sum is to be applied
does not overlap any component of the link.  Choose $n$ crossings of
$C$ with itself such that flipping the parities of these $n$ crossings
(i.e. over to under or vice versa) makes $C$ a $0$-framed unknot (this
step relies on the winding number condition).  Apply move II $n$ times to
create a Hopf link for each of these crossings, then apply move I$'$
twice at each crossing as in Figure \ref{fg:crossingswitch}(A).  The effect
of these moves is to make $C$ a $0$-framed unknot (of course the disk
it bounds intersects $L$ in many places), and in this new link, the
band connected sum of $B$ with $C$ along the same band is an instance of
move I$'$.
After applying this band connected sum, the vicinity of each
of the $n$ crossings looks like the left-hand side of Figure
\ref{fg:crossingswitch}(B), and undoing the two instances of move I$'$
and the instance of move II at each crossing results in the right side
of Figure \ref{fg:crossingswitch}(B), which is a projection of the original band
connect sum of $B$ with $C.$
\end{proof}

\begin{rem}
Surgery is described by an unoriented framed link while the quantities
of the upcoming sections
will naturally be invariants of oriented framed links.  Thus our
strategy will be to find an oriented  link invariant which is
unchanged by reversal of the orientation of any component, as well as
by spin Kirby moves I$'$ and II.  Theorem \ref{th:Kirby}  and
Proposition \ref{pr:spin} then say that such an invariant will
actually be an invariant of the framed three-manifold presented by
surgery on that link.
\end{rem}
\section{Reshetikhin-Turaev Type Spin Invariants}

\subsection{Modularity and spin modularity}

We review here the main results of \cite{\nonsimply}  which will be
relevant to the 
question of spin invariants. 

\begin{figure}[hbt]
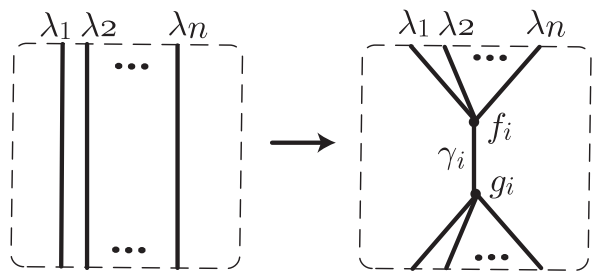

$$\pic{binding}{100}{-45}{0}{0}$$
\caption{Binding edges together}\label{fg:bind}
\end{figure}

Recall (see for example Kirillov \cite{Kirillov96}, Turaev
\cite{Turaev94} and  \cite{Sawin96a,\nonsimply}) that the quantum group
$U_q(\g),$ where $\g$ is a simple Lie 
algebra and $q$ is a root of unity, forms a ribbon Hopf algebra.  More 
precisely, 
the set of representations of $U_q(\g)$ spanned by the finite collection of
irreducible representations with highest weight in the Weyl alcove
forms a semisimple ribbon $*$-category with the truncated tensor
product $\trunc.$  Thus we get  a numerical invariant of framed graphs  with edges 
labeled by such representations and vertices labeled by invariant
elements of an appropriate tensor product of the labels of incident
edges and the duals of those labels (invariant elements are in general
 morphisms in the ribbon category)  such that the following hold:
\begin{enumerate}

\item the invariant of a graph with an edge labeled by $\lambda \oplus
\gamma$ is the sum of the invariants of the same graph with that edge
labeled by $\lambda$ and $\gamma$ respectively, the labels on the
adjacent vertices being projected appropriately,
\item if the label of an edge is replaced by an isomorphic object and
the labels of the adjacent vertices are composed with the isomorphism
in the obvious way, the invariant is unchanged.   In particular, link
components can be unambiguously labeled by elements of $\Gamma,$ rather
than objects,
\item  the invariant of a graph with an edge labeled by the trivial
object (the weight $0$) is
the same as the invariant of the graph with that edge  deleted,
\item  the
invariant of a graph with an edge labeled by $\lambda$ is the
invariant of the graph with the orientation of that edge reversed
and the label replaced by the dual $\lambda^\dagger,$ the labels of the adjacent
vertices remaining the same,
\item  the invariant of a graph with
a link  component labeled by $\lambda \trunc \gamma$ is the invariant of
the graph with that component replaced by two parallel components
(according to the framing) labeled by $\lambda$ and $\gamma$
respectively,
\item the invariant of the connected sum of two graphs 
along  edges labeled by a simple object $\lambda$ 
 is the product
of the invariants of the two graphs divided by $\qdim(\lambda)$  and
\item for any objects $\lambda_1, \ldots, \lambda_n$ there is a
collection of pairs of invariant elements $f_i \in \lambda_1 \trunc
\cdots \trunc \lambda_n \trunc \gamma_i$ and $g_i \in
\lambda_1^\dagger \trunc \cdots \lambda_n^\dagger \trunc
\gamma_i^\dagger$ for various simple objects $\gamma_i$ such that if $L$
is any graph and a ball intersects $L$ as in Figure \ref{fg:bind}, the
sum of the invariants of the graphs on the right side of the figure
equals the invariant of $L.$
\end{enumerate} 

The invariant of the unknot labeled by a representation $\lambda$ is
called the quantum dimension of $\lambda,$ $\qdim(\lambda),$ and
the invariant of the Hopf link labeled by $\lambda$ and $\gamma$ is
called $S_{\lambda,\gamma}.$ $C_\lambda$ is the modulus one complex
number by which  the link invariant is multiplied  when a component
labeled by $\lambda$ is given a positive full twist.

An object $\lambda$ (i.e. a representation) in the ribbon category is
called \emph{degenerate} if 
$R_{\lambda,\gamma}= R_{\gamma,\lambda}^{-1}$ for all objects
$\gamma,$ where $R_{\lambda,\gamma}$ is the morphism corresponding to
a  positive crossing.  By necessity if $\lambda$ is degenerate
then $C_\lambda=\pm 1.$  The set of degenerate objects with
$C_\lambda=1$ (such objects are called \emph{even}) forms a symmetric
subcategory which is isomorphic to the 
representation theory of some compact group, and M\"uger proves
\cite{Muger99} that one can quotient by any full subcategory of this
subcategory in the sense that there is a minimal semisimple ribbon $*$-category
which admits a ribbon $*$-functor from the original category to it
sending all the objects in the subcategory and only those objects to direct sums of the  trivial object.
In particular, if we quotient by the subcategory of \emph{all}
even objects  the resulting category will have no even degenerate
objects except sums of the trivial object.

Suppose $\Cat$ is a semisimple ribbon $*$-category with the property
that the only even degenerate simple object is  trivial, and
suppose that it contains a degenerate simple object $\mu$ with
$C_\mu=-1$ (naturally, we call $\mu$ an \emph{odd} degenerate object).   Then
$\mu \trunc \mu$ is degenerate and even, and therefore is a sum of copies of the
trivial object.   Since $\mu$ is simple, this means $\mu^\dagger=\mu$ and
$\mu \trunc \mu$ is exactly the trivial object (here $\mu^\dagger$ is
the dual object to $\mu$).  If $\nu$ is a
degenerate object with $C_\nu=-1,$ then $\mu \trunc \nu$ is even, so it
is a sum of trivial objects, so $\nu$ is a sum of objects isomorphic
to $\mu.$  Thus we have the following.

\begin{thm}
Let $\Cat$ be a semisimple ribbon $*$-category, and let $\Cat'$ be the
quotient by the full symmetric subcategory of even degenerate objects
as described above.  Then there are two possibilities
\begin{alist}
\item If all degenerate objects in $\Cat$ are even, then $\Cat'$
contains no degenerate objects, and therefore according to M\"uger is
modular, and can be used to construct a three-manifold invariant and
TQFT following  Reshetikhin and Turaev \cite{RT91,Turaev94}.
\item If $\Cat$ contains any degenerate objects which are not even, then
$\Cat'$ contains exactly one simple degenerate object $\mu,$ and $\mu$
satisfies $C_\mu=-1$ and $\mu \trunc \mu$ is trivial.  In this case we
call $\Cat'$ \emph{spin modular,} and will construct spin and 
compatibly framed 
three-manifold invariants from it.
\end{alist}
\end{thm}
\begin{rem}  The $*$-structure on the ribbon category is 
not necessary.   Brugui\`eres has a similar construction to M\"uger's
which replaces the $*$-structure with semisimplicity together with a fairly easy
to check and general condition on the degenerate objects.   The
$*$-structure, which is available in all the cases of interest to us
and clearly is related to the physical origin of the relevant
invariants, is merely a convenience in this situation.
\end{rem}

\subsection{Framed and spin three-manifold invariants}

Let $\Cat$ be a spin modular category, and let $\mu$ be the simple
degenerate object.

\begin{lem}
For each simple object $\lambda \in \Cat,$ the object 
$\sigma(\lambda)\defequals \mu \trunc
\lambda$ is  a simple object distinct from $\lambda,$ and in any link
with an even framed component labeled by $\lambda,$ the invariant is
unchanged if that label is replaced by $\sigma( \lambda).$  In
particular, if $V$ is the vector space of formal linear combinations
of simple objects in $\Cat$ (the link invariant extends by linearity to links
labeled by elements of $V$) then $\sigma$ extends to a $\Z_2$ action
on $V,$ the link invariant descends to an 
invariant of even links labeled by elements of $V/\sigma,$ the map $V \to
V/\sigma$ preserves the duality map $\dagger$ (thought of as another
$\Z_2$ action  on $V$) and finally the $S$-matrix gives a nondegenerate
pairing on $V/\sigma.$
\end{lem}

\begin{proof}
Since $\mu \trunc (\mu \trunc \lambda)=\lambda$ is simple, $\mu \trunc
\lambda$ must be simple.  Of course $C_{\sigma(\lambda)}=-C_\lambda$
is different from $C_\lambda$ so $\lambda$ and $\sigma(\lambda)$ are
distinct.  For any  link $L,$ the invariant of $L$ with an even 
component labeled by $\mu \trunc \lambda$ is the invariant of $L$ with
that component doubled and labeled by $\mu$ and $\lambda$
respectively.   Since $\mu$ is degenerate we have
$R_{\mu,\gamma}=R_{\gamma,\mu}^{-1},$ so in particular in any
projection of this link any crossing involving $\mu$ can have its
parity switched  (over to under or vice versa) without changing the
invariant.  A sequence of such changes can unlink and unknot the
component labeled by $\mu,$ so the invariant of the doubled link is
equal to the invariant of the original link time $\qdim(\mu)C_\mu^n,$
where $n$ is the framing of the component labeled by $\mu.$  Now $n$
is even and $C_\mu=-1,$ so we get only the factor of $\qdim(\mu).$
But since $\mu \trunc \mu$ is trivial, we know $\qdim(\mu)^2=1,$ and
since in a $*$-category quantum dimensions are positive,
$\qdim(\mu)=1.$

Of course $\sigma$ extends to a $\Z_2$ action  on $V$ by linearity, which
commutes with the duality map $\dagger.$  
Again by linearity it is true that labeling a component by any $v \in
V$ gives the same value to the invariant of an even link as labeling
it by $\sigma(\lambda),$ so we can as well label components of an even
link by equivalence classes $\{v,\sigma(v)\} \in V/\sigma.$  Thus the
pairing defined by
$\bracket{\lambda,\gamma}=S_{\lambda,\gamma},$ descends to a well-defined pairing on $V/\sigma.$

Now the tensor product $\trunc$ extends by linearity to an
associative, distributive multiplication with identity on $V$ (here we
identify $\lambda \oplus \gamma$ with the vector  $\lambda + \gamma$
of $V$), thus
making $V$ an algebra.   Since
$$\sigma(\lambda)\trunc \gamma = \sigma(\lambda \trunc \gamma)=\lambda
\trunc \sigma(\gamma),$$
 the quotient $V/\sigma$ inherits the algebra
structure.

For each equivalence class
$[\lambda]$ of a simple object $\lambda$ notice
$f_{[\lambda]}([\gamma])\defequals \bracket{\lambda,\gamma}/\qdim(\lambda)$
is a nontrivial homomorphism from $V/\sigma$ to $\C.$  Supposing the
pairing $\bracket{\, ,}$  is degenerate,
then these $\dim(V)/2$ homomorphisms must be linearly dependent, and
thus two must be equal.  If $f_{[\lambda]}=f_{[\lambda']},$ then
$\bracket{\lambda,\gamma}/\qdim(\lambda) =\bracket{\lambda',\gamma}/\qdim(\lambda')$ for all
$\gamma.$  Now since $\bracket{\lambda' \trunc
\lambda^\dagger,\gamma}=\bracket{\lambda',\gamma}
\bracket{\lambda^\dagger,\gamma}/\qdim(\gamma)$ we have
$$\bracket{\lambda\trunc\lambda^\dagger,\gamma}/\qdim(\lambda) =
\bracket{\lambda'\trunc\lambda^\dagger,\gamma}/\qdim(\lambda')$$
 for all $\gamma.$  In
particular, since $f_{[\id]}$ (where $\id$ is the trivial object,
which is the multiplicative identity) is
nontrivial there is a minimal 
idempotent $\omega$ such that $f_{[\id]}(\omega)=1$ but
$f_{[\gamma]}(\omega)=0$ if $f_{[\gamma]} \neq f_{[\id]}.$  Now
$$\bracket{\lambda \trunc \lambda^\dagger,\omega}=\sum_\gamma
N_{\lambda,\lambda^\dagger}^\gamma
\bracket{\gamma,\omega} \geq N_{\lambda,\lambda^\dagger}^\id =1,$$
 where
$N_{\lambda,\gamma}^\delta$ is the multiplicity of $\delta$ in
$\lambda \trunc \gamma.$  Thus 
$$\bracket{\lambda'\trunc
\lambda^\dagger,\omega}\geq \qdim(\lambda')/\qdim(\lambda)>0$$
 so there
exists a $\gamma$ such that 
$N_{\lambda',\lambda^\dagger}^\gamma =1$ and $f_{[\gamma]}=f_{[\id]}.$  

Of course if $[\gamma]=[\id]$ then $\lambda'=(\lambda^\dagger)^\dagger=\lambda$
or $\lambda'=\sigma(\lambda^\dagger)^\dagger=\sigma(\lambda),$ so since $[\lambda']
\neq [\lambda]$ we conclude that there is a $\gamma \neq \id,\mu$ such
that $\bracket{\gamma,\lambda}=\qdim(\gamma) \qdim(\lambda)$ for 
all $\lambda.$  M\"uger
proves that this property implies $\gamma$ is degenerate, so we reach a
contradiction and conclude the pairing was nondegenerate.
\end{proof}

\begin{lem}\label{lm:omega} Let $\omega=\sum_\gamma \qdim(\gamma) \gamma.$
\begin{alist}
\item For all $v \in V,$  $v \trunc \omega=\qdim(v)\omega$
\item  $\bracket{\lambda,\omega}/\qdim(\lambda)$ 
 is $\qdim(\omega)$ if $\lambda=\id$ or
$\lambda=\mu,$ and $0$ otherwise.
\item $\bracket{\omega,\omega}=2\qdim(\omega)\neq 0.$
\end{alist}
\end{lem}

\begin{proof}\mbox{}
\begin{alist} 
\item It suffices to prove this for $v=\lambda$ with  $\lambda \in \Lambda.$
\begin{multline*}
\lambda\trunc  \omega = \sum_\gamma \qdim(\gamma) \lambda \trunc \gamma= \sum_{\gamma,\eta}
\qdim(\gamma) N_{\lambda,\gamma}^\eta \eta\\
=\sum_{\gamma,\eta} N_{\lambda,\eta^\dagger}^{\gamma^\dagger} \qdim(\gamma)
\eta =\sum_{\gamma,\eta} N_{\lambda,\eta^\dagger}^{\gamma^\dagger}
\qdim(\gamma^\dagger)\eta\\
=\sum_\eta \qdim(\lambda \trunc \eta^\dagger) \eta=\qdim(\lambda) \sum_\eta
\qdim(\eta) \eta=\qdim(\lambda) \omega.
\end{multline*}
\item By the previous point, for any simple $\gamma$
$$
\bracket{\lambda,\omega}=\bracket{\lambda,\omega \trunc
\gamma}/\qdim(\gamma)\\
=\bracket{\lambda,\gamma} \bracket{\lambda,\omega}/(\qdim(\gamma)
\qdim(\lambda))
$$
so either $\bracket{\lambda,\omega}=0$ or for every $\gamma$ we have
$\bracket{\lambda,\gamma}=\qdim(\lambda) \qdim(\gamma).$  The second
condition we have already noted is equivalent to the degeneracy of
$\lambda,$ so this only happens when $\lambda=\id$ or $\lambda=\mu.$ 
In both cases the formula follows immediately.
\item Using the previous point
\begin{multline*}
\bracket{\omega,\omega}=\sum_\gamma \qdim(\gamma)
\bracket{\gamma,\omega}\\
=\qdim(\id)\bracket{\id,\omega}+\qdim(\mu)\bracket{\mu,\omega}=2\qdim(\omega).
\end{multline*}
Of course $\qdim(\omega)=\sum_\gamma \qdim(\gamma)^2>0.$
\end{alist}
\end{proof}

\begin{thm} \label{th:spininvt}
Let $L$ be an even link representing a three-manifold $M$ with a compatible
framing and let $I(L)$ be the invariant associated to the spin-modular
category $\Cat$ acting on  $L$ with each component
labeled by $\omega/\sqrt{2\qdim(\omega)}.$  I.e., if $L$ has components
$1,\ldots,n$ and $L_{\lambda_1,\ldots,\lambda_n}$ is $L$ with the $n$
components labeled by simple objects $\lambda_1,\ldots,\lambda_n$
respectively and $F$ is the link invariant, then 
\begin{equation} \label{eq:tminvt}I(L)=\left(\frac{1}{\sqrt{2\qdim(\omega)}}\right)^n
\sum_{\lambda_1,\ldots,\lambda_n} \left(\prod_{i=1}^n
\qdim(\lambda_i)\right) F(L_{\lambda_1,\ldots,\lambda_n})
\end{equation}
where the sum is over isomorphism classes of simple objects and
$\qdim(\omega)=\sum_\lambda \qdim(\lambda)^2.$  Then $I(L)$ 
is invariant under the spin Kirby moves and therefore is an
invariant of $M$ and its framing.
\end{thm}

\begin{proof}
Notice first that 
$$\omega^\dagger = \sum_\gamma \qdim(\gamma) \gamma^\dagger=
\sum_\gamma \qdim(\gamma^\dagger) \gamma^\dagger=  \omega,$$
so $I(L)$ is unchanged by reversing the orientation of any component
of $L.$

For invariance under Kirby move II as pictured in Figure
\ref{fg:Kirby}, notice by Properties 3 and 6  of the link 
invariant the value of $F$ on the link on the right is
$\bracket{\omega,\omega}/2\qdim(\omega)=1$ times that on the left,
using Lemma \ref{lm:omega}(c).

For invariance under Kirby move I$'$ as pictured in Figure
\ref{fg:Kirby}, the argument is given pictorially in Figure
\ref{fg:invariance}, where the first equality is by Property 7 of the
invariant, the second by Lemma
\ref{lm:omega}(b), the third by the degeneracy of $0$ and $\mu,$ the
fourth by Lemma \ref{lm:omega}(b) and the fifth by Property 7 again.
\end{proof}

\begin{figure}[hbt] 
\begin{multline*}
F\left(\pic{kirby1}{65}{-30}{-2}{2}\right) = \sum_i
 F\left(\pic{kirby2}{65}{-30}{-2}{2}\right)=\\ 
\sum_{i:\gamma_i=\id,\mu}
F\left(\pic{kirby3}{65}{-30}{-2}{2}
\right)\qdim(\omega) 
= \sum_{i:\gamma_i=\id,\mu}
F\left(\pic{kirby4}{65}{-30}{-2}{2}\right)\qdim(\omega)\\
=\sum_i F\left( \pic{kirby5}{65}{-30}{-2}{2} \right) =
F\left(\pic{kirby6}{65}{-30}{-2}{2}\right) 
\end{multline*}
\caption{Invariance under Kirby move I$'$}\label{fg:invariance}
\end{figure}

\begin{prop}
If $\Cat$ is a ribbon $*$-category whose even  degenerate objects form
a cyclic group of invertibles and $F$ is the associated link invariant then
Equation (\ref{eq:tminvt}) applied to the category $\Cat$
is an invariant of compatibly framed three-manifolds as above.
\end{prop}

\begin{proof}
The proof of Proposition 2 of \cite{\nonsimply} shows that if $\Cat$ is
any ribbon $*$-category and $\Cat'$ is the quotient as in
\cite{Bruguieres99,Muger99} of $\Cat$ by a cyclic group of invertible
even degenerate objects, then the image of the expression $I(L)$ above
under the functor from $\Cat$ to $\Cat'$ is the same expression in the
image category (the statement of the proposition discusses  only the case when the
quotient is by the full set of degenerate objects, but the argument
does not use this fact in any way).  Thus if the full set of even
degenerate objects is a 
cyclic group of invertible elements and $\Cat'$ is therefore modular, the
formula $I(L)$ gives an invariant of $2$-framed three-manifolds as in
\cite{\nonsimply}, and if not then $\Cat'$ is spin-modular and $I(L)$
gives an invariant of compatibly framed three-manifolds as in the previous theorem.
\end{proof}

\begin{rem}
Lemma \ref{lm:omega}(a) applied to a particular knot gives the
invariance of $\omega$ under band connect sum of the unknot with that
knot, and with a little more effort of any knot separated by a sphere with
that knot.  It does not appear to give the full invariance under move
I  (which would be conceptually superior to our indirect proof
through move I$'$) without substantially more effort.   This effort
would amount to switching focus from the space associated to the torus
(roughly what we call $V$) to the space associated to the \emph{punctured}
torus.  The appropriate setting for this would be the full axioms of
an extended TQFT of Walker \cite{Walker??},  whose generalization to
the spin category offers a very interesting question.
\end{rem}
\begin{cor}
Every closed subset $\Lambda$ of the Weyl alcove, for every quantum group
$U_q(\g)$ at every level $k$ (i.e., for $q$ an arbitrary root of
unity) except possibly the exceptions listed below, yields a
framed (or $2$-framed) three-manifold invariant by the formula $I(L)$ 
above, where $F$ is the standard quantum group link invariant.  The 
possible exceptions are:  $\g$ having Dynkin diagram $D_{2n}$ and the set
$\Gamma_{\Z_2 \times \Z_2}$ of weights in the root lattice, where the 
group of degenerates is not cyclic,  as well as 
at level $k=2$ the exceptional sets for $D_n$ and $B_n$  
 discussed in \cite{\closed}, where the group associated to the 
 subcategory of  degenerates is not 
 commutative.   
\end{cor}

\begin{proof}
By \cite{\closed}, every closed subset of the Weyl alcove (closed
meaning that the truncated tensor product of two elements of the set
is the sum of elements of the set) yields a semisimple ribbon $*$-subcategory of the
standard ribbon category associated to the Weyl alcove whose 
degenerate objects are invertible and  by \cite{\nonsimply} form (except 
for cases mentioned) a cyclic group.  In \cite{\nonsimply} and \cite{\closed} it is
determined when these closed subsets include odd degenerate objects,
and thus whether the formula $I(L)$ gives an invariant of $2$-framed
three-manifolds or compatibly framed three-manifolds.
\end{proof}

\begin{rem}
    Presumably Formula (\ref{eq:tminvt}) is preserved by the 
    quotient, and thus defines a $2$-framed or compatibly framed 
    three-manifold invariant, even when the subcategory of even 
    degenerate objects is not generated by a cyclic group, but a proof 
    is currently lacking.
\end{rem}

\subsection{Decomposition into prime invariants}

In \cite{\nonsimply} it was found that many of the closed subsets of the
Weyl alcove which give modular quotients and hence $2$-framed
three-manifold invariants could be factored as a product of 
others from the list.  In this subsection we analyze the
decomposition  of subsets whose quotients are spin-modular.

Let $\g$ be  a simple Lie algebra and  $G$ be the simply-connected compact
Lie group with Lie algebra $\g,$ and   consider the Weyl alcove at
level $k.$  If   $Z$ is a subgroup of the center
$Z(G)$ of $G,$   let $\Gamma_Z$ be the intersection of the Weyl
alcove with those weights in the Weyl chamber which are the highest
weights of representations of $G$ on which $Z$ acts trivially.  Let
$\Delta_Z$ be the  set of weights  which are $k$ times a
fundamental weight $\lambda$ whose inner product with every element of
$\Gamma_Z$ is an integer.  In fact, there is a bijection $\ell$ from
$Z(G)$ to a certain subset of the fundamental weights,  and $\Delta_Z$ is the image $k 
\ell[Z]$ of $Z,$ where $k$ represents the linear map multiplication by
$k$ in the lattice.
Then  \cite{\nonsimply} shows that the sets 
$\Gamma_Z$ and $\Delta_Z$ are closed. 

Recall from \cite{\nonsimply} that if $\Gamma$ is the set of isomorphism
classes of simple objects of a ribbon category $\Cat$ (in particular
if it is a closed subset of the Weyl alcove) then we say $\Gamma$ is
the product of two subsets $\Gamma'$ and $\Gamma''$ if 
\begin{enumerate}
\item the intersection $\Gamma' \cap \Gamma''$ consists of even
degenerate objects,
\item the product $\trunc$ of any element of $\Gamma'$ with an element
of $\Gamma''$ is simple (i.e. is an element of $\Gamma$),
\item every element of $\Gamma$ is a product of an element of
$\Gamma'$ and $\Gamma''$ and
\item if $\lambda' \in \Gamma'$ and $\lambda'' \in \Gamma''$ then
$C_{\lambda' \trunc \lambda''}= C_{\lambda'} C_{\lambda''}.$
\end{enumerate}

\begin{prop}
If $\Gamma$ is the product of $\Gamma'$ and $\Gamma''$ then
\begin{equation} \label{eq:product}
I_\Gamma(L)=I_{\Gamma'}(L)I_{\Gamma''}(L)
\end{equation}
where I is the invariant of the previous subsection computed in the
categories associated to $\Gamma,$ $\Gamma',$ and $\Gamma''$
respectively.
\end{prop}

\begin{proof}
As in the previous subsection let $V_\Gamma$ be the formal vector
space spanned by isomorphism classes of simple objects in $\Cat,$ and
let $V_{\Gamma'}$ and $V_{\Gamma''}$ be the corresponding vector
spaces for $\Cat'$ and $\Cat''.$  Of course the truncated tensor
product gives an algebra homomorphism $\phi:V_{\Gamma'} \tensor
V_{\Gamma''} \to V_\Gamma$ which by point 3 in the definition above is onto.
 It is shown in \cite{\nonsimply} that the link invariant with a
component labeled by $\phi(a \tensor b)$ is the product of the link
invariants with components labeled by $a$ and $b$ respectively.  From
this we see that $\phi(\omega' \tensor \omega'')$ has the property
that $\phi(\omega' \tensor \omega'')v=\phi(\omega' \tensor
\omega'')\qdim(v),$ because $\omega' \tensor \omega''$ has this
property in $V_{\Gamma'} \tensor
V_{\Gamma''}.$
Thus
\begin{multline*}\phi(\omega' \tensor \omega'')\qdim(\omega) 
= \phi(\omega' \tensor
\omega'')\omega \\
=  \qdim(\phi(\omega' \tensor \omega'')) \omega=
\qdim(\omega') \qdim(\omega'') \omega,
\end{multline*}
the last equality being a consequence of the behavior of the invariant
under $\phi.$  Thus $\phi(\omega' \tensor \omega'')$ and $\omega$ are nonzero
multiples of each other, since $\qdim(\omega),$ $\qdim(\omega'),$ and
$\qdim(\omega'')$  are all nonzero.  Since  $I(L)$ is easily seen to
be unchanged
if we replace $\omega$ by a nonzero multiple, the result follows from
the invariant's behavior under $\phi.$ 
\end{proof}

The following is proven in \cite{\nonsimply}.

\begin{prop}\label{pr:product}
Suppose $Z \subset Z'$ are subgroups of the center of $G,$
 $\Gamma'$ is the closed subset
generated by $\Gamma_{Z'}$ and $\Delta_Z,$ and
$Z_0=Z \cap (k \circ \ell)^{-1} [\Gamma_{Z'}].$ Then $\Gamma' \subset
\Gamma_{Z_0}$ is of the form $\Gamma'=\Gamma_Y$ for some $Y$ and
$\Delta_{Z_0}=\Delta_Z \cap \Gamma_{Z'}$ consists of 
degenerate invertible objects for $\Gamma'.$  If all of $\Delta_{Z_0}$ is even then 
$\Gamma'$ is the  product of $\Gamma_{Z'}$ and $\Delta_Z.$ These are
the only cases in which $\Gamma$ decomposes
into a product, apart from $D_{2n}.$  \qed
\end{prop}

We shall be particularly interested in theories coming from closed
subsets $\Gamma_Z$ with the property  that $\Delta_Z \subset \Gamma_Z.$  The
following proposition shows that all other closed sets $\Gamma_Z$ appear as factors
of these, and the next proposition gives conditions for when a closed
set has this property.

\begin{prop}\label{pr:DW}
If $\g\neq D_{2n},$ every $\Gamma_{Z'}$ is a factor of a  $\Gamma_{Z_0}$ with the
property that $\Delta_{Z_0} \subset \Gamma_{Z_0}.$ 
\end{prop}

\begin{proof}
Given $\Gamma_{Z'},$ let $Z$ be $Z'$ if all $\lambda$ in  $\Delta_{Z'} \cap \Gamma_{Z'}$
satisfy $C_\lambda=1.$  If $\Delta_{Z'} \cap \Gamma_{Z'}$ contains any
elements with $C_\lambda=-1,$ then those with $C_\lambda=1$ form an
index two subgroup, and since $\Delta_{Z'}$ is cyclic, $\Delta_{Z'}$ too must
have an index two subgroup whose intersection with $\Gamma_{Z'}$ has
only objects with $C_\lambda=1.$  In that case let $Z$ be the corresponding index
two subgroup of $Z'.$  Then $\Delta_{Z'}$ is generated by $\Delta_Z$
and any element of $\Delta_{Z'} \cap \Gamma_{Z'}$ with $C_\lambda=-1.$
Thus in particular $\Delta_{Z'}$ is generated by $\Delta_Z$ and 
elements of $\Gamma_{Z'}.$  So whatever $Z$ is $\Gamma_{Z'}$ and
$\Delta_Z$ generate the same set $\Gamma'$ generated by $\Gamma_{Z'}$
and $\Delta_{Z'}.$  The previous proposition shows $\Gamma_{Z'}$
and $\Delta_Z$ generate $\Gamma_{Z_0}$ and that $\Delta_{Z_0}=\Delta_{Z}
\cap \Gamma_{Z'}$ consists of even degenerate objects for
$\Gamma_{Z_0},$ so $\Gamma_{Z_0}$ is the product of $\Gamma_{Z'}$ and
$\Delta_Z.$ 

\end{proof}

\begin{prop}\label{pr:DWtype}
 $\Gamma_Z$ contains $\Delta_Z$ as degenerate objects if and only if
\begin{alist}
\item  $k(\ell(z),\ell(z))$ is an even integer for all $z \in Z,$ in which case 
$\Gamma_Z$ contains only even degenerate objects or
\item  $k(\ell(z),\ell(z))$ is an integer for all
$z\in Z,$ with at least one of these integers odd, in which case $\Gamma_Z$
contains an odd degenerate object and  $Z$ contains
an index $2$ subgroup.
\end{alist}
\end{prop}

\begin{proof}
Recall from \cite{\nonsimply} $\lambda \in \Gamma_Z$ if and only if
$(\lambda,\ell(z))
\in
\Z$ for all
$z \in Z.$  since $Z$ is cyclic, this is equivalent to $(\lambda,\ell(z)) \in
\Z$ for $z$ a generator.  So $\Delta_Z \in \Gamma_Z$ is equivalent to
$k(\ell(z'),\ell(z)) \in \Z$ for all $z' \in Z$ and for some generator $z\in
Z,$ which is to say $k(\ell(z),\ell(z))\in \Z$ for some generator $z \in Z,$
or equivalently for all $z \in Z.$  

Now it is shown in \cite{\nonsimply} that 
$$R_{\lambda,k\ell(z)}R_{k\ell(z),\lambda}^{-1}= e^{2\pi i
k(\ell(z),\lambda)},$$
so $k\ell(z)$ with $z \in Z$ is always degenerate for $\Gamma_Z$ if it is in
$\Gamma_Z.$  thus we need only check the additional assertions in (a) and (b).
 Again from \cite{\nonsimply}, $C_{k\ell(z)}=\exp(\pi i
k(\ell(z),\ell(z))),$ so this is one if and only if $k(\ell(z),\ell(z))$ is
even.  Finally, notice the set of $z \in Z$ such that $k(\ell(z),\ell(z))$ is
even forms  a subgroup.  If it is proper it has index two, because the
product of two elements not in this subgroup is in this subgroup.
\end{proof}

Finally, we observe that when $\Delta_Z$ is cyclic and contains an odd degenerate
object,  the framed invariant $I_{\Delta_Z}$ that results is actually an
ordinary
$2$-framed invariant.  In fact it is the invariant Murakami, Ohtsuki and Okada
associate to the quotient of  $Z$  by the entire group of degenerates \cite{MOO92}.
Specifically,  let $\lambda$ be a generator of $\Delta_Z,$ and suppose $N$ is the
least such that $\lambda^N$ is a degenerate object. Recalling the
above formulas for $C_{k\ell(z)}$ and $R_{k\ell (z),\gamma}
R^{-1}_{\gamma,k\ell(z)}= C_{k\ell(z) \trunc \gamma}C_{k\ell(z)}^{-1}
C_\gamma^{-1},$   to  say that 
$\lambda^n$ is degenerate is to say that
$C_{\lambda^{n+1}}=C_{\lambda^n}C_\lambda,$ since it suffices to check
the degeneracy condition against a generator.  Now $\lambda=k\ell(z)$
 and $\lambda^n=k\ell(z^n)$ for some $z\in Z,$ so if we let 
 $r=\exp(\pi i (k\ell(z), \ell(z)))$ then $C_\lambda=r$ and 
 $$R_{\lambda^n,\lambda} R_{\lambda,\lambda^n}= e^{2\pi i k 
 (\ell(z),\ell(z^n))}=r^{2n}.$$
 So by induction $C_{\lambda^n}=\lambda^{n^2}.$  Thus  $N$ 
 is the least such that $\lambda^N$ is degenerate if and only if  $N$ is 
 the least  $N$ such that  $r^{2N}=1.$ 

If $\Delta_Z$ contains an odd degenerate object then $\lambda^N$ is
odd, and thus $C_{\lambda^N}=r^{N^2}=-1.$  We conclude
that $N$ is odd and $r$ is a primitive $2N$th root of unity.  

 Now quotient by the even degenerate objects so that $Z$
becomes $\Z_{2N}$ and we can 
identify the simple object $\lambda^m$ with the number $m \in \Z_{2N}.$
We claim that the invariant of a link with $n$
components labeled by the entries in $\vec{l}=(l_1,\ldots,l_n) \in
(\Z_{2N})^n$  is $r^{\vec{l}^tA \vec{l}},$ where $A$ is the linking
matrix of the link.  To see this notice the formula agrees with the
invariant for a link composed of $n$ unlinked but possibly framed unknots, and that the
invariant and the proposed formula both change by  $r^{2pq}$ when a
crossing between components labeled by $p$ and $q$ respectively
switches parity.

Thus the compatibly framed  invariant  is
$$I_{\Delta_Z}(L) = (2\sqrt{N})^{-n} \sum_{\vec{l}
\in (\Z_{2N})^n} r^{\vec{l}^tA \vec{l}}.$$
But if $A$ is an even matrix then $\vec{l}^tA \vec{l}$ is always even,
so we may replace $r$ by the primitive $N$th root of unity $-r$
without  changing the value.  In that case the contribution from any
$\vec{l}$ is the same as from $\vec{l} + (0,0,\ldots, 0,N,0,\ldots,0)$
and thus we may take the sum over $\Z_N$ at the cost of multiplying by
$2^n,$ and thus get
$$I_{\Delta_Z}(L) = (\sqrt{N})^{-n} \sum_{\vec{l}
\in (\Z_N)^n} (-r)^{\vec{l}^tA \vec{l}}.$$
Now notice this is exactly the invariant Murakami, Ohtsuki and Okada
call $Z(-r,N),$ which is a $2$-framed three-manifold invariant.   
Thus  it does not depend on the spin structure.

\subsection{Relationship to geometry and physics}

In \cite{DW90}, Dijkgraaf and Witten discuss under what circumstances
one  expects a TQFT and three-manifold invariant to arise from the
Chern-Simons field theory of a (possibly nonsimply-connected) compact 
simple Lie
group $G.$  Their approach is to 
define the Chern-Simons functional when the principal $G$-bundle over
the three-manifold is not trivial in terms of a cohomology class in
$H^4(BG,\Z).$  The action is computed by choosing a four-manifold and 
principal bundle (or more generally a homology class in $H_4(BG,\Z)$)
bounded by the three-manifold and principal bundle and pairing the
cohomology class with the fundamental class of the four-manifold.
This result must always be an integer when the three-manifold is
trivial in order for the path integral to be
well-defined, which is why $H^4(BG)$ must be  
taken with integer coefficients and the reason for the integrality
conditions on $k$ derived by Dijkgraaf and Witten.

Dijkgraaf and Witten add an intriguing point.  If the group is such
that the generating class of $H^4(BG,\Z)$ (Recall $H^4(BG,\Z)\cong \Z$ for 
every  compact simple
group $G$ except the one associated to the Dynkin diagram $D_{2n}$ with
trivial center) when integrated against the fundamental class of a
\emph{spin} four-manifold is even, then half-integer multiples of the
generator (which we may think of as $H^4(BG)$ classes with
half-integer coefficients) would give a well-defined action if the 
four-manifold is forced  to be spin.  Thus in these cases we expect
half-integer cohomology classes to give spin TQFTs (as in the ordinary 
case the theory still depends on the signature of the extending 
four-manifold, i.e., the $2$-framing, so really we expect compatibly 
framed TQFTs).  Dijkgraaf and
Witten show that $\mathrm{SO}(3)$ in particular has this property and discuss
the expected spin TQFT in this case  (which should occur when the
$\mathrm{SU}(2)$ level is $2$ modulo $4$).

In fact the property is quite general and occurs exactly in the
situation covered by Proposition \ref{pr:DWtype}.
\begin{prop}
Suppose $G$ is a compact Lie group such that $H^4(BG,\Z)=\Z$ and
$\pi_1(G)$ has an index two subgroup (i.e. $g$ admits a double cover).  Then the
generating class of $H^4(BG,\Z)$ pulled back via a $G$-bundle on a
spin four-manifold has even integral.
\end{prop}

\begin{proof}
Since $H^4(BG,\Z)=\Z$ it follows $H^4(BG,\Z_2)=\Z_2.$  We claim there
is a $w \in H^2(BG,\Z_2)$ such that $w \wedge w$ is nontrivial and
hence is equal to the unique nontrivial element of $H^4(BG,\Z/2).$
Specifically, $w$ is the obstruction to lifting a $G$ bundle to a
$\tilde{G}$ bundle, $\tilde{G}$ being its double cover, which
necessarily exists because of the condition on $\pi_1(G).$

To see that $w \wedge w$ is nontrivial, consider $\C\mathrm{P}^2,$ the
four-manifold  constructed by the surgery process of Section 2 from
the $+1$-framed unknot.  More precisely, it is formed by attaching a
single $2$-handle to a $0$-handle along the framed unknot and then
attaching a $4$-handle to the resulting $S^3$ boundary.  We construct a
$G$-bundle over $\C\mathrm{P}^2$ as follows.  Attach the trivial bundle
over the $2$-handle to the trivial bundle over the $0$-handle via an
overlap map on the boundary which is homotopic to an
element of $\pi_1(G)$ not in the index two subgroup.
  Extend this bundle over the $4$-handle.  Of course this
bundle does not lift to $\tilde{G},$ so the image $\tilde{w}$ of $w$
in $H^2(\C\mathrm{P}^2,\Z_2)=\Z_2$ is nontrivial.  The image of $w \wedge
w$ is $\tilde{w} \wedge \tilde{w}$ which is nontrivial because the
intersection pairing  given by the $1 \times 1$ identity linking
matrix is nondegenerate.

On
four-manifolds with even forms  (which includes all spin
four-manifolds) $\int_M w \wedge w=0 \bmod 2,$ so the integral of 
every class of $H^4(BG,\Z)$ against such four-manifolds is even.  
\end{proof}

\begin{cor}
The techniques of Dijkgraaf and Witten predict a spin (really a 
compatibly framed) Chern-Simons
theory  at the levels and groups given in Proposition \ref{pr:DWtype}(b).
\end{cor}

\begin{proof}
Dijkgraaf and Witten's techniques, in conjunction with the previous
proposition, predict a spin TQFT associated to the Lie group $G/Z$
with half-integer levels, 
where $G$ is a simply-connected group and $Z$ is a subgroup of the
center containing an index two subgroup.     Now
integer levels of $G/Z$ correspond as levels of $G$ to integer
multiples of $N,$  $N$ being the least integer such that
$N(\ell(z),\ell(z))$ is an even integer for all $z \in Z.$  So
half-integer levels correspond to odd multiples of $N/2.$  

Thus Dijkgraaf and Witten predict a spin TQFT for $G/Z$ at level $k$ if
$k(\ell(z),\ell(z))$ is an integer for all $z\in Z$ (because it is a
half-integer multiple of an even number) but at least one of these
numbers is odd (otherwise $N(\ell(z),\ell(z))$ is a multiple of $4$
for all $z \in Z,$ and $N$ would not be the least such meeting the defining
condition).

\end{proof}

From Dijkgraaf and Witten's point of view we should interpret these spin
Chern-Simons theory as theories we compute on spin four-ma\-ni\-folds, and
perhaps even as invariants of spin four-manifolds which happen to depend
only on their boundary (and signature).   By defining our invariants
in terms of surgery, we have partially modeled this feature.  By
assigning numbers to even links, we are really assigning numbers to
certain spin four-manifolds, and discovering the fact that the numbers
depend only on the spin boundary and signatures.  Of course our
numbers are only assigned to spin four-manifolds which admit a handle
decomposition as one $0$-handle and some $2$-handles.  It is an
interesting question whether from a spin-modular category one can
naturally associate to \emph{every}  spin four-manifold an invariant
which reduces to this one when the four-manifold admits a handle
decomposition as described.  Presumably this would involve assigning
some sort of a label to link components representing $1$-handles
(Kirby's dotted circles \cite{Kirby89}).   This question may shed light on
efforts to construct interesting four-manifold invariants from
algebraic structures related to quantum groups.

\subsection{Identification with spin invariants of Kirby and Melvin}

A special case of the construction of the preceding subsections is the
group $\mathrm{SO}(3)$ which is $\mathrm{SU}(2)/\Z_2.$  Here the center
contains a single nontrivial element $z$ and $\ell(z)$ is the unique
fundamental weight with $(\ell(z),\ell(z))=1/2$ and thus $N=4.$  So we
expect spin $\mathrm{SO}(3)$ invariants at $\mathrm{SU}(2)$ levels which are $2$
modulo $4.$  These are exactly the levels at which Kirby and Melvin
\cite{KM91} construct
invariants of spin three-manifolds from the representation theory of
quantum $\mathrm{su}_2$  (see also Blanchet \cite{Blanchet92}, as well as
\cite{BM96}, where Blanchet and 
Masbaum define a spin TQFT giving this invariant).   We will show that
Kirby and Melvin's invariant is exactly the $\mathrm{ SO}(3)$ invariant, and in
fact that every spin invariant of the previous subsection can be
computed in a manner analogous to theirs using the subset of
representations associated to the double cover.

Recall that Kirby and Melvin present a framed three-manifold by a link $L$ with
a \emph{characteristic sublink} $C$ such that for each component $a$ of $L$
$$\sum_{c \in C} a\cdot c=a\cdot a \bmod 2$$
where the dot represents the linking number between the two components
(or in the case of $a\cdot a,$ the framing).  This corresponds to the
three-manifold obtained by surgery on the link $L,$ with the unique
spin structure which when restricted to $S^3-L,$  extends over the
components of $L-C$ but not over the components of $C.$   We describe
a generalization of their invariant and show that it gives exactly 
the spin invariants we constructed in the previous subsections.

Let $\Delta_Z \subset \Gamma_Z$ be such that $\Gamma_Z$ contains odd
degenerate objects.  Then as in Proposition \ref{pr:DWtype}, the
subset $Z_0= \{z \in Z: k(\ell(z),\ell(z)) \in 2\Z\}$ 
is an index two subgroup, and therefore $\Gamma=\Gamma_{Z_0}$
has  $\Gamma_Z$ as a subset and contains only even degenerate objects.

Let 
\begin{align*}
\omega&= \sum_{\gamma \in \Gamma} \qdim(\gamma) \gamma\\
\omega_0&= \sum_{\gamma \in \Gamma_Z} \qdim(\gamma) \gamma\\
\omega_1&= \omega-\omega_0.
\end{align*}

\begin{prop}\label{pr:KMinvariant}
If $L$ is a link with characteristic sublink $C,$ let $F(L,C)$ be the
invariant (in $\Gamma$) of $L$ with every component of $C$ labeled by
$\omega_1$ and every other component of $L$ labeled by $\omega_0,$ let
$U_+$ and $U_-$ be the invariants of the respectively $+1$  and $-1$
 framed unknots labeled by $\omega_1,$ and let $n$ be the number
of components of $L.$
Then
\begin{equation} \label{eq:KMframed}
J(L,C)=F(L,C)/(U_+ U_-)^{n/2}
\end{equation}
is an invariant of the framed (i.e. $2$-framed and spin)  manifold determined by $(L,C)$
and
\begin{equation} \label{eq:KMspin}
J'(L,C)=\left( \frac{U_-}{U_+}\right)^{\sigma/2}J(L,C)
\end{equation}
is an invariant of the ordinary spin manifold determined by $(L,C).$
\end{prop}

\begin{proof}
We shall confirm invariance of the second quantity, that of the first
follows.  Invariance under orientation reversal is clear.  

According to Kirby and Melvin an invariant of a link with a
characteristic sublink is an invariant of the spin manifolds if it is
unchanged by the following move:  add a $\pm 1$ framed unknot to the
link (possibly linking with other components),  apply a positive or negative
full twist to the  disk it bounds  (so as to change the linking matrix
of the link) and add it to the characteristic sublink if and only if
the sum of its linking numbers with the existing characteristic
sublink  is even. 

Notice first that as argued earlier the formula $J'(L,C)$ applied to a link
gives the same answer as the same formula interpreted in the quotient
of $\Gamma$ by $Z_0,$ so we may assume that $\Gamma$ is modular and
that $\Gamma_Z$ is spin modular.  Let us call the unique odd
degenerate object in $\Gamma_Z$  $\mu.$  

Observe first that if $\lambda$ is a simple object of $\Gamma
-\Gamma_Z,$ then $S_{\lambda,\mu}=-\qdim(\lambda).$  To see this, 
note that since $\mu^2$ is the trivial object and
$\qdim(\lambda)=S_{\lambda,\mu \trunc \mu}=
S_{\lambda,\mu}^2/\qdim(\lambda)$ we have that
$S_{\lambda,\mu}=\pm \qdim(\lambda).$  But supposing $S_{\gamma,\mu}=\qdim(\lambda),$ then
if $\gamma$ is an other simple object in $\Gamma$ but not in $\Gamma_Z,$
then $\lambda \trunc \gamma \in \Gamma_Z,$ so that $\qdim(\lambda)
\qdim(\gamma)= S_{\lambda \trunc \gamma,\mu}=S_{\lambda,\mu}
S_{\gamma,\mu} =\qdim(\lambda) S_{\gamma,\mu}$ and we conclude
$S_{\gamma,\mu}=\qdim(\gamma).$  But this is certainly true for
$\gamma \in \Gamma_Z,$ so $\mu$ is degenerate for
$\Gamma.$   This is a contradiction so $S_{\lambda,\mu}=-\qdim(\lambda).$
 
Consider the result of Kirby and Melvin's move, and suppose first
that the unknot is to be added to the characteristic sublink, because
its linking number with the old sublink is even.   The invariant
$F(L,C)$ is a sum over labelings of the components of $L,$ with those
in $C$ labeled by elements of $\Gamma-\Gamma_Z$ and those not in $C$ labeled
by elements of $\Gamma_Z.$  Choose such a labeling, and consider the invariant of
this labeled link with a particular label $\kappa$ on the new unknot.   The
condition on the linking number means that the new unknot surrounds a
collection of strands an even number of which have labels not in
$\Gamma_Z,$ and therefore the tensor product of all these labels is a
sum of labels in $\Gamma_Z$ (here we use Property 7 of the ribbon 
category).  If $\kappa$ is in $\Gamma_Z,$ then
$\kappa \trunc \mu$ is a distinct label with $C_{\kappa \trunc
\mu}=-C_\kappa$ and $S_{\gamma,\kappa}=S_{\gamma,\kappa \trunc
\mu}$ for $\gamma \in \Gamma_Z$ and we see that labeling the new
unknot by $\kappa$ versus $\kappa \trunc \mu$ contributes the same
amount with opposite sign to the computation of the total invariant.  Thus
labeling the new unknot by $\omega_0$ would give a total invariant of
$0,$ so that in the computation of $J'(L,C)$ we might as well replace
the label $\omega_1$ on the new unknot with $\omega.$  The same
argument applies to $U_+$ and $U_-,$ so we see that in this case the
invariance of $J'(L,C)$ under this move is equivalent to the
invariance of the ordinary manifold invariant of $\Gamma$ under this
move.

Now suppose that the new unknot is not to be added to the
characteristic link, because its linking number with the
characteristic link is odd.  The same argument as above shows the
truncated tensor product of the labels of the strands going through
the new unknot is a sum of elements of $\Gamma-\Gamma_Z.$  If
$\kappa$ is a label for the new unknot which is not in $\Gamma_Z$ and
we compare the effect on the invariant of replacing $\kappa$ by
$\kappa \trunc \mu,$ we see that $C_{\kappa\trunc
\mu}=C_{\kappa},$ but now an unknot labeled by $\mu$ surrounding
a sum of labels not in $\Gamma_Z$ contributes $-1,$ so again $\kappa$
and $\kappa \trunc \mu$ contribute opposite amounts to the sum (in
this case they may not be distinct, but then  $\kappa$ contributes zero)  and thus
labeling the unknot by $\omega_1$ results in a total invariant of
zero.   Once again we may as well replace the label of $\omega_0$
with $\omega,$ and again the results follows from the invariance of
the standard invariant under the move.
\end{proof}

In the case $\g=\mathrm{ su}_2,$ $k \equiv 2 \mod 4,$ and $Z=\Z_2,$ we
have that $\Gamma$ is the whole set of representations, $\Gamma_Z$ is
the set of integer spin representations, and our formula reduces
exactly to Kirby and Melvin's formula  (they sum only over half the
representations, but using their symmetry principle,  this is
equivalent to summing over all the representations, as they note).

That this invariant is the one we already constructed is now obvious
by taking the characteristic sublink to be empty.

\begin{prop} \label{pr:KM=me}
The invariant $J(L,C)$ when applied to a compatibly framed
three-manifold, gives the same result as the invariant $I$ associated to $\Gamma_Z.$
\end{prop}

\begin{proof}
If we present the compatibly framed manifold by an even link $L,$ then
notice that the empty link is a characteristic sublink, and
$(L,\emptyset)$ is the Kirby-Melvin presentation of this $2$-framed
spin three-manifold.  Of course $F(L,\emptyset)=F(L)$ as defined in
the second subsection of this section, and thus the invariants are
equal as long as the normalizations are equal, that is if $U_+
U_-=F(H),$ with $H$ the Hopf link labeled by $\omega_0.$  But of
course the Hopf link represents $S^3$ with the standard spin structure
and $2$-framing, so $J(H,\emptyset)=1,$ which means $F(H)=U_+ U_-.$
\end{proof}

\begin{rem}
Kirby and Melvin's argument that the sum of the $\Gamma_Z$ invariant
over all possible spin structures on a given manifold adds up to the
$\Gamma$ invariant of the manifold goes through in the general case by
an analogous argument.
\end{rem}

\begin{rem}
In the Kirby-Melvin formulation of the invariant, any $2$-framed spin
three-manifold can be represented by a link and characteristic
sublink, and thus we get an invariant of $2$-framed spin-manifolds,
without the compatibility constraint.  This extension of our  framed invariant
amounts to a canonical choice of a sixteenth root of the invariant of
the Kummer surface, and thus the additional information it detects is
at most Rohlin's invariant.   Whether there is a canonical way to do 
this in a
more general ribbon $*$-category (i.e., a category that does not 
already come embedded with index $2$ in a category without odd 
degenerate objects), or in the situation of the next
section, is an open question.
\end{rem}

\begin{rem}
If we consider $\mathrm{su}_2$ at level $k=2,$ then
$\Delta_{\Z_2}=\Gamma_{\Z_2}$ is a closed subset consisting of the
trivial object and an odd degenerate object.   The invariant
$I_{\Delta_{\Z_2}}$ is the Murakami et al's invariant $Z(1,1),$ which
is completely trivial  (it assigns $1$ to all $2$-framed manifolds).
However, the Kirby-Melvin extension is nontrivial, as it depends on
Rohlin's invariant (in fact it is $\exp(-3 \pi i \mu(M)/8)$).  
\end{rem}

\section{Hennings Type Spin invariants}

\subsection{Even link invariants from quasitriangular Hopf algebras}

Recall that a \emph{quasitriangular Hopf algebra} is a Hopf algebra 
$$(H,\cdot,
1, \Delta, \epsilon, S)$$
(so that $H$ is an algebra with
multiplication $\cdot$ and unit $1,$  the dual space $H^\dagger$ is an
algebra with multiplication $\Delta^\dagger$ and unit $\epsilon,$ 
$\Delta\negmedspace : \negmedspace H
\to H \tensor H$ is
a homomorphism, and $S\negmedspace :\negmedspace H \to H$ is an antihomomorphism satisfying $x^{(1)}
S(x^{(2)})=S(x^{(1)})
\tensor x^{(2)} = \epsilon(x) 1,$
where $\Delta(x)=x^{(1)} \tensor x^{(2)}$ is Sweedler's index-saving
notation) together with an element $R=\sum_i a_i \tensor b_i \in H
\otimes H$ such that
\begin{equation} \label{eq:R} 
\Delta^\mathrm{op}(x) R = \sum_i (x^{(2)} a_i \tensor x^{(1)} 
b_i)= \sum_i (a_i x^{(1)} \tensor b_i x^{(2)}) =R \Delta(x)
\end{equation}
\begin{multline} \label{eq:lcable}
(\Delta \tensor \id_H)(R)= \sum_i(a_i^{(1)} \tensor a_i^{(2)} \tensor
b_i) \\
= \sum_{i,j} (a_i \tensor a_j \tensor b_ib_j) = R_{1,2}
R_{2,3}
\end{multline}
\begin{multline} \label{eq:rcable}
(\id_H \tensor \Delta)(R)=\sum_i a_i \tensor b_i^{(1)} \tensor
b_i^{(2)} \\= \sum_{i,j} (a_i a_j \tensor b_j \tensor b_i)=R_{1,3} R_{1,2}.
\end{multline}

From these relations  follow all sorts of useful facts, including
\begin{multline} \label{eq:YB}
R_{1,2}R_{1,3} R_{2,3}= \sum_{i,j,k} (a_ia_j \tensor b_i a_k \tensor
b_jb_k) \\
= \sum_{i,j,k}(a_j a_k \tensor a_ib_k \tensor b_i b_j)=
R_{2,3} R_{1,3} R_{1,2}
\end{multline}
\begin{equation} \label{eq:trivial}
(\epsilon \tensor \id_H)(R)=\sum_i \epsilon(a_i)b_i = 
1 =
\sum_i \epsilon(b_i) a_i =(\id_H \tensor \epsilon)(R)
\end{equation}
\begin{multline} \label{eq:inverse}
(S \tensor \id_H)(R) = \sum_i S(a_i) \tensor b_i = \\
R^{-1} = \sum_i
a_i \tensor S^{-1}(b_i) = (\id_H \tensor S^{-1})(R)
\end{multline}
\begin{equation}
(S \tensor S)(R)= \sum_i S(a_i) \tensor S(b_i) = \sum_i a_i \tensor
b_i =R.
\end{equation}

If $H$ also possesses an element $g \in H$ such that $\Delta g = g \tensor
g,$ $\epsilon(g)=1,$ $g^{-1} u g^{-1}=S(u)$ and $S^2(x)=gxg^{-1}$ for
all $x \in H,$ where $u=\sum_i S(b_i) a_i,$ then $H$ is called a
\emph{ribbon Hopf algebra.}   Its representation
theory forms a ribbon category, and associated to it by the recipe of
 Ohtsuki \cite{Ohtsuki??} or
Kauffman and Radford \cite{KR95b} is a numerical 
invariant of framed links with  components labeled by \emph{quantum
characters.}  Here a quantum character of $H$ is a functional $f \in H^\dagger$ such
that $f(\Ad_x(y))=\epsilon(x) f(y),$ where $\Ad_x(y)=x^{(1)} y S(x^{(2)}).$
  Equivalently (according to Drinfel'd \cite{Drinfeld90}), a quantum
character is  a functional $f$ such that $f(xy)=f(yS^2(x))$ for all $x,y \in H.$

Here we  associate to a finite-dimensional quasitriangular (but not 
necessarily ribbon) Hopf
algebra $H$ an invariant of \emph{even} links with components labeled by
quantum characters.  The invariant is defined precisely in analogy
with the definition for ribbon categories (in fact, if the
quasitriangular Hopf algebra is extended in the usual way of
Reshetikhin and Turaev \cite{RT90} to a ribbon
Hopf algebra, the invariants agree on even links), and we will imitate
the construction of \cite{KR95b} closely.

\begin{figure}[hbt]
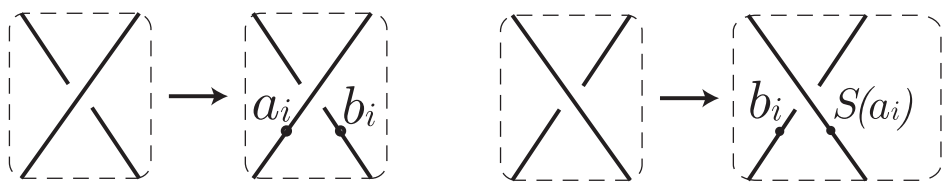

$$\pic{crossing}{65}{-30}{0}{0}$$
\caption{Labels assigned to crossings} \label{fg:crossing}
\end{figure}

Let $L$ be an even  link with a quantum character $\lambda_j$ associated to
each component $C_j$ of $L.$  Choose a presentation of $L$ with every
component having winding number one, and also with a height function.
We associate
to the projection a collection of \emph{decorated projections,} with Hopf
algebra elements assigned to various noncritical points on each
component as follows.  For each crossing, choose one value of the
index $i$ in $R=\sum_i a_i \tensor b_i$ (we call such a choice for all
crossings a \emph{state}), and label points near the crossing as shown
in Figure \ref{fg:crossing}, where the vertical in the figure is the height function.  Now choose a base point on each
component which is not a crossing or a critical point of the height function. Let us
call the direction along the component which at the basepoint is of
increasing height function the \emph{basepoint direction.}
To each noncritical point of each component, associate
an \emph{integer rotation number} by counting all the critical points
visited when traveling from the base point to that point in the
basepoint direction, counting
each critical point which rotates clockwise as $-1$ and each which
rotates counterclockwise as $+1.$  That is, the nearest integer to the
actual rotation from the basepoint to that point divided by $\pi.$  

For each state, assign a number to each component of the  associated decorated link
projection as follows. Travel from the base point
once around the component in the basepoint direction.  Form the product of all the decorated
Hopf algebra elements on the component, left to right in the order of visitation,
each acted on by $S$ raised to the integer rotation number.  The
number assigned to the component is  $\lambda_i$ of this product
if the basepoint direction coincides with  the orientation of the
component and is $\lambda_i$ of $S$ of the product if they are in
opposite directions.  The number assigned to the state is the product
of the numbers assigned to the components, and finally the number
assigned to the projection is the sum of over all states of the
numbers assigned to the states.

 Given a state, or any decorated link diagram, it is possible to compute the
number assigned to that state `in stages.'  Specifically, a
\emph{fragment} of a link projection is the intersection of the
projection with a disk whose boundary intersects the projection
transversely and not at double points or critical points.  Choose for
each strand of the fragment a basepoint (several strands may belong
to the same component of the entire link).  Now  each decorated
point in the fragment has associated to  it an integer rotation
relative to the basepoint of its strand.  Assign to each strand the
product of $S^{n_i}(x_i),$ where $x_i$ is a decoration, $n_i$ is the
integer rotation number, and the product is over all decorations.  The
order of the product comes from traversing the strand in the basepoint
direction, traveling 
from one endpoint to the other for open strands, and from the
basepoint once around for closed strands.

The key observation is that if we erase all decorations in the
fragment and replace them with the chosen basepoints decorated with the
assigned products as above, the number associated to the entire state
is unchanged.  

In particular, notice that if we define the set of states of a fragment to be
the set of all possible assignments of index values to the crossings
in the fragment, then the set of states of the projection is the
product of the set of states of the fragment and the set of states of
its complement (which is also a fragment).   So the number associated
to the entire projection, which is a sum over states, can be found
by summing over all pairs of states a number computed by evaluating
the state on the fragment and on its complement and combining appropriately.

What's more,  if a fragment of a projection is replaced with a
perhaps topologically distinct fragment, but such that 
the strands connect the same pairs of boundary points, the rotation
numbers of the boundary points from the basepoint are the same, and
the sum of the decorations of the basepoints, viewed as an element of
$H^{\tensor n}$ where $n$ is the number of open strands, is the same,
then the projection with the replaced fragment will yield 
the same number as the original projection.  The proof of the
following proposition offers a concrete illustration of this
observation, where it is the key to proving invariance under the
regular isotopy moves.

\begin{prop}\label{pr:invariance}
The number computed above depends neither on the specific projection
nor on the choice of basepoints, but only on the even link and quantum
characters. We call this invariant quantity $K(L),$  or $K_{\lambda_1,
\ldots, \lambda_n}(L)$ when the labels require explicit mention.
\end{prop}

\begin{proof}
For invariance under choice of basepoint, it suffices to check the
quantity is unchanged if one basepoint is moved past a crossing or a
critical point. This fact is generally true for decorated projections
and does not rely on special properties of the $R$-matrix.

To see that it is unchanged when a basepoint is moved past a crossing
or any decorating Hopf algebra element,  consider for a
given state and component, the product $PS^{\pm 2}(x),$ where $x$ is the
last decoration encountered in the traversal about the component and
$P$ is the product of the rest of 
the decorating elements with appropriate powers of $S.$  The $\pm 2$
represents the integer rotation number of the point decorated by $x,$
which because of the winding number condition on the projection is $2$
if the direction of traversal agrees with the orientation and $-2$ if
it disagrees.  In either case the defining property of the quantum
character gives
\begin{align*}
\lambda_i(PS^2(x))&=\lambda_i(xP),\\
\lambda_i(S(PS^{-2}(x)))&=\lambda_i(S^{-1}(x)S(P))=\lambda_i(S(P)S(x))=\lambda_i(S(xP)),
\end{align*}
so the computed quantity is the same as for the case where the base
point is just below instead of just above the decorated point.

To see that it is unchanged when the base point is moved past a
critical point, consider a base point at which the orientation of the
component is  upwards, which is to be moved past a minimum of
the height function  configured so that the orientation rotates
clockwise around it  (the other cases all work similarly). Viewing the
rest of the link as a fragment with the fragment basepoint just above
the basepoint for the component in the  link  and decorated by $x,$ the
value of the decorated link projection before moving the basepoint is $\lambda_i(x),$ and after
moving it is $\lambda_i(S(S^{-1}(x)))=\lambda_i(x).$

To see that the invariant is independent of the choice of projection,
notice by Proposition \ref{pr:Reidemeister} we may check that it does
not change under  regular isotopy and the height function
moves.

The first regular isotopy move reduces to the equation 
$$\sum_{i,j} S(a_i) a_j \tensor b_i b_j=1 \tensor 1 =\sum_{i,j} a_i S(a_j)
\tensor b_i b_j$$
 which is a restatement of  Equation (\ref{eq:inverse}).  The second isotopy
move is exactly the Yang-Baxter equation (\ref{eq:YB}).  The height
function moves are immediate from the definition.  
\end{proof}

\begin{rem}
In the presence of a ribbon structure the trace with respect to a
representation $V$ of the Hopf algebra can be made into a quantum
character by adding the charmed element, $\lambda_V(\cdot)=\tr(g^{-1}
\cdot).$  The link invariant described in this section would then
correspond to the usual Reshetikhin and Turaev link invariant
associated to $V$ \cite{RT90}.  For a nonsemisimple Hopf
algebra there are typically other quantum characters.   But even in
the absence of the ribbon structure, for every representation $V$ we
can always choose an intertwiner from $V^\tensor V^\dagger$ to the
trivial representation (if $V$ is irreducible, it is unique up to
scale) and viewing elements of the Hopf algebra via the representation
as elements of $V \tensor V^\dagger,$ we get a functional on the Hopf 
algebra which
proves to be a quantum character.  Thus even in the merely 
quasitriangular situation
representations give (even) link invariants.  We expect
that there is a categorical structure, analogous to but weaker than
the notion of ribbon category, which axiomatizes the structure that
allows the category of representations of such a Hopf algebra to give
even link invariants.   It seems plausible that the quotient of a
ribbon category containing odd degenerate objects by the full set of
degenerate objects, which does not make sense as a ribbon category,
would be well-defined as a category of this sort.   We conjecture
further that associated to the half-integer level Chern-Simons
theories of the previous section there are quasitriangular but not
ribbon Hopf algebras such that an appropriate truncation of their
representation theory  (corresponding to the truncated tensor product
construction for ordinary quantum groups) yields the spin
Chern-Simons theories at those levels. 
\end{rem}

A few important facts about $K$  will be used
in the next section, all following easily from the form of the computation:
\begin{enumerate}
\item  $K(L_1 \cup L_2)=K(L_1)K(L_2),$ where $L_1 \cup L_2,$ the
distant union of the labeled links $L_1$ and $L_2,$ is the link formed by
embedding each into $S^3$ so that they are separated by a sphere.
\item Reversing the orientation of a component corresponds to
composing the quantum character labeling that component with $S.$
\end{enumerate}

\subsection{Integrals and the three-manifold invariant}

This section relies heavily on some general results about
finite-dimensional Hopf algebra.  A reference that contains everything
we need is Radford \cite{Radford94}.

Recall that a left (respectively right) integral in the dual of  a Hopf algebra $H$ is a
functional $\lambda \in H^\dagger$ such that $\gamma \lambda=
\gamma(1)\lambda$ (respectively $\lambda
\gamma=\gamma(1) \lambda$) for all $\gamma \in H^\dagger$.   We
say $H$ is unimodular if there is a $\lambda$ which is a left and
right integral simultaneously, in which case it is unique up to scale.
If $\lambda$ is a left and right integral, then by \cite[Theorem 3]{Radford94}, $\lambda(ab)=\lambda(bS^2(a))$ for all
$a,b \in H,$ so that $\lambda$ is a quantum character.   Also by
Proposition 3 of the same article $\lambda \circ S=\lambda,$ so that
(assuming now that $H$ is quasitriangular and thus determines a link
invariant as described in the previous subsection)
changing the orientation of a component labeled by $\lambda$ does not
change the invariant.   

The  integral $\lambda$ enjoys a particularly distinctive property: It
generates $H^\dagger$ as a free $H$-module.  More specifically, if we
define for each $h \in H$ the functional $f_h \in H^\dagger$ by
$$f_h(a)=\lambda(ah)$$
for all $a \in H,$ then the map $h \mapsto f_h$ is one-to-one and
onto.  From this it follows that $f_\Lambda=\epsilon,$ where $\Lambda$
is an appropriately normalized left (and therefore right) integral for
$H.$  

Now consider, for $H$ quasitriangular with $R=\sum_i a_i \tensor b_i,$
a map from $H^\dagger$ to $H$ given by
$$D:f \mapsto \sum_{i,j}f(a_i b_j)b_ia_j.$$
Notice by Equation (\ref{eq:R}) this is a homomorphism.  We say $H$
is \emph{factorizable} if this
homomorphism $D$ (the Drinfel'd map) is bijective. In that case the
image of the integral $\lambda$ must be a nonzero integral for $H.$
Since $\lambda(D(\lambda)) =f_{D(\lambda)}(1)\neq 0,$ we can normalize $\lambda$ once
and for all by the condition $\lambda(D(\lambda))=1.$  

\begin{thm}\label{th:Hennings}
If $H$ is a unimodular, factorizable quasitriangular Hopf algebra the
invariant $K(L)$ of a link $L$ all of whose components are labeled  by
the quantum character $\lambda$ is an invariant of the compatibly framed
three-manifold described by surgery on $L.$
\end{thm}

\begin{figure}[hbt]
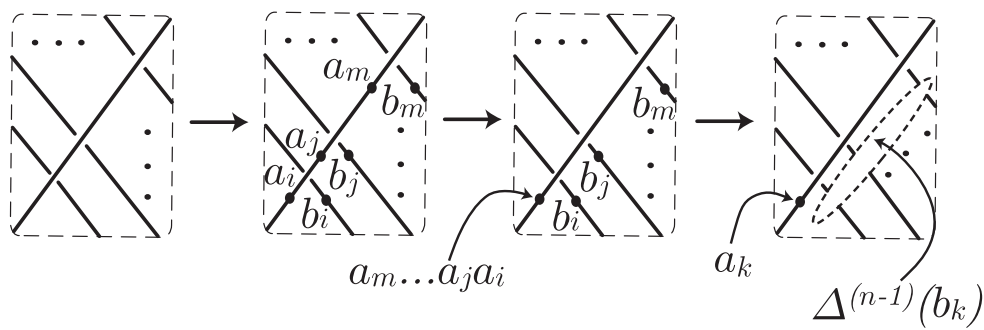

$$\pic{cable}{120}{-50}{0}{0}$$
\caption{Hennings invariant of a cabled crossing}\label{fg:cable}
\end{figure}

\begin{proof}
We have already argued that this invariant is unchanged by orientation
reversal of any component.  Invariance under spin Kirby move II
follows from the normalization of $\lambda$ and Property 1 of the
list at the end of the previous subsection (the invariant of the Hopf 
link is $\lambda((D\lambda))$).

Before proving invariance under move I$'$, consider the fragment in
Figure \ref{fg:cable}.   By repeated application of Equation
 (\ref{eq:rcable})  we see that the sum over all states of the
values assigned to the fragment is $\sum_k a_k \tensor
\Delta^{n-1}(b_k),$ where the $n$ entries of the tensor product
label the components from left to right at the bottom of the fragment.
Similarly, using Equations (\ref{eq:lcable}) and (\ref{eq:rcable}), we
see that the same fragment with the opposite parity is associated to
$\sum_i S(b_i) \tensor \Delta^{n-1}(a_i),$ the mirror image of this
fragment is associated to $\sum_i \Delta^{n-1}(S(b_i)) \tensor a_i,$ and
the mirror image with opposite parity by $\sum_i \Delta^{n-1}(a_i)
\tensor b_i.$  

With this in hand we see that the fragment on the left-hand side of
Figure \ref{fg:spinFR} corresponds  (with the same convention of
decorating the strands from left to right at the bottom of the
picture) to
\begin{multline*}
\sum_{i,j,k}  a_i \tensor \lambda(a_jb_k) \Delta^{n-1}(b_ka_jb_i)
=\sum_i a_i \tensor \Delta^{n-1}(D(\lambda)b_i)\\
=\sum_i a_i \tensor \epsilon(b_i)\Delta^{n-1} (D(\lambda))
=1 \tensor \Delta^{n-1}(D(\lambda))
\end{multline*}
the last line by Equation (\ref{eq:trivial}).   Of course, replacing
each $a_i$ in the above computation by $S(b_i)$ and every $b_i$ with
$a_i,$ we see that the right-hand side equals the same thing, and thus
the two fragments are interchangeable in the calculation of the invariant.
\end{proof}

The most important special case of this construction arises when we
consider the Drinfel'd double $H_A$ of a finite-dimensional Hopf
algebra $A.$  It is well-known \cite{Radford94} that $H_A$ is
quasitriangular, unimodular and factorizable.  Thus to every
finite-dimensional Hopf algebra $A$ this construction associates an
invariant of compatibly framed  (or spin) three-manifolds.  There is
ample empirical evidence, but as yet no proof, that this is exactly the
invariant associated by Kuperberg to the Hopf algebra $A$ in
\cite{Kuperberg96}.


\bibliographystyle{plain}

\end{document}